\documentclass[twoside]{article}
\usepackage[accepted]{aistats2015}

\usepackage{hyperref,url,breakurl,fancybox,amsmath,amssymb,amsthm,bm,paralist,booktabs,subfig,graphicx}%fancybox => rounded box, breakurl => url-s can break
\usepackage[mathscr]{eucal}
\usepackage{empheq}

\newcommand{\X}{\mathscr{X}}   %domain of the distributions
\renewcommand{\H}{\mathscr{H}} %RKHS(K)
\newcommand{\Hmeta}{\mathscr{G}} %RKHS-meta
\newcommand{\Kmeta}{k_{\Hmeta}} %RKHS-meta
\newcommand{\R}{\mathbb{R}}    %real number
\newcommand{\Z}{\mathbb{Z}}%   integer
\newcommand{\E}{\mathbb{E}}    %expectation
\renewcommand{\Pr}{\mathbb{P}} %probability
\renewcommand{\L}{\mathscr{L}} %linear operator

\newcommand{\Eo}{\mathcal{E}}  %E-objective
\renewcommand{\d}{\mathrm{d}}  %d(x): integration
\renewcommand{\b}{\mathbf}     %bold

\newcommand{\A}{\mathscr{A}}
\newcommand{\B}{\mathscr{B}}
\newcommand{\N}{\mathscr{N}}
\newcommand{\M}{\mathscr{M}}%probability measure

\renewcommand{\P}{\mathcal{P}}%P(b,c)
\newcommand{\tb}{\textbf}%bold text
\newcommand{\Bo}{\mathcal{B}}%Borel sets
\renewcommand{\S}{\mathscr{S}}%sigma-alg

\DeclareMathOperator*{\argmin}{arg\,min}

\newtheorem*{theorem}{Main theorem (bound on the excess risk)}
\newtheorem{consequence}{Consequence}
\newtheorem{lemma}{Lemma}

\newtheorem{definition}{Definition}

\setdefaultleftmargin{0em}{0em}{0em}{0em}{0em}{0em}%own
\begin{document}

\runningtitle{Two-stage sampled learning theory on distributions}
\runningauthor{Szab{\'o} et al.}

\twocolumn[
\aistatstitle{Two-stage sampled learning theory on distributions}
\aistatsauthor{\hspace{0.25cm}Zolt{\'a}n Szab{\'o}$^1$ \hspace*{1cm} Arthur Gretton$^1$ \hspace*{1cm} Barnab{\'a}s P{\'o}czos$^2$ \hspace*{1cm} Bharath Sriperumbudur$^3$}
\aistatsaddress{\hspace*{0.2cm} $^1$Gatsby Unit, UCL \hspace*{1.3cm} $^2$Machine Learning Department, CMU \hspace*{1.3cm} $^3$Department of Statistics, PSU}]

\begin{abstract}
We focus on the distribution regression problem: regressing to a real-valued response from a probability distribution.
Although there exist a large number of similarity measures between distributions, very little is known about their
generalization performance in specific \emph{learning tasks}. Learning problems formulated on distributions have an inherent
two-stage sampled difficulty: in practice only  samples from sampled distributions are observable, and one has to build an estimate on similarities computed between sets of points.
To the best of our knowledge, the only existing method with consistency guarantees for distribution regression 
requires kernel density estimation as an intermediate step (which suffers from slow convergence issues in high dimensions), and the domain of the distributions to be
compact Euclidean. In this paper, we provide theoretical guarantees for  a remarkably simple algorithmic alternative to solve the distribution regression problem: embed the distributions to a reproducing kernel Hilbert space, and learn a ridge regressor from the embeddings to the outputs.
Our main contribution is to prove the consistency of this technique in the two-stage sampled setting under mild conditions (on separable, topological domains endowed with kernels). 
For a given total number of observations, we derive convergence rates as an explicit function of the problem difficulty. As a special case, we answer a $15$-year-old open question: we establish 
the consistency of the classical set kernel [Haussler, 1999; G{\"a}rtner et.\ al, 2002] in regression, and cover more recent kernels on distributions, 
including those due to [Christmann and Steinwart, 2010].

\end{abstract}

\vspace*{-0.2cm}
\section{INTRODUCTION}
\vspace*{-0.1cm}
We address the learning problem of \emph{distribution regression} in the two-stage sampled setting \cite{poczos13distribution}: we  regress from probability measures to real-valued responses, where we only have bags of samples from the probability distributions.
Many classical problems in machine learning and statistics can be analysed in this framework. On the machine learning side, multiple instance learning \cite{dietterich97solving,ray01multiple,dooly02multiple} can be thought of in this way,
in the case where each instance in a labeled bag is an i.i.d.\ (independent identically distributed) sample from a distribution. On the statistical side, tasks might include point estimation of statistics on a distribution (e.g., its entropy or a hyperparameter),
where a supervised learning method can help in parameter estimation problems without closed form analytical expressions, or if simulation-based results are computationally expensive.

Before reviewing the existing techniques in the literature, let us start with a somewhat informal definition of the distribution regression problem, and an intuitive phrasing of 
our goal. Let us suppose that our data consist of $\b{z}=\{(x_i,y_i)\}_{i=1}^l$, where $x_i$ is a probability distribution, $y_i \in \R$, and
each $(x_i,y_i)$  pair is i.i.d.\ sampled from a meta distribution $\mathcal{M}$. However, we do not observe $x_i$ directly; rather, we observe a sample
$x_{i,1},\ldots, x_{i,N_i} \stackrel{i.i.d.} {\sim} x_i$. Thus the observed data are $\hat{\b{z}} = \{(\{x_{i,n}\}_{n=1}^{N_i},y_i)\}_{i=1}^l$. 
Our goal is to predict a new $y_{l+1}$ from a new batch of samples $x_{l+1,1},\ldots, x_{l+1,N_{l+1}}$ drawn from a new distribution $x_{l+1}$. 
For example, in a medical application the $i^{th}$ patient might be identified with a probability distribution ($x_i$), which can be periodicly accessed, measured by 
blood tests ($\{x_{i,n}\}_{n=1}^{N_i}$). We are also given some health indicator of the patient ($y_i$), which might be inferred from his/her blood measurements.
Based on the observations ($\hat{\b{z}}$), one might try to learn the mapping from the set of blood tests to the health indicator; and the hope is that by 
observing more patients (larger $l$) and performing a larger number of tests (larger $N_i$) the estimated mapping ($\hat{f}=\hat{f}({\hat{\b{z}}})$) becomes more ``precise''.

The performance of the estimated mapping ($\hat{f}$) depends on the assumed function class ($\H$), the family of $\hat{f}$ candidates. Let 
$f_{\H}$ denote the best estimator from $\H$ given infinite training samples ($l=\infty$, $N_i=\infty$), and let $\Eo[f_{\H}]$ be its prediction error.
Our goal is to obtain upper bounds for the $0\le \Eo[\hat{f}]-\Eo[f_{\H}]$ quantity which hold with high probability. More precisely, we are aiming at
\begin{compactenum}
  \item deriving upper bounds on the excess risk, proving consistency: We construct $\Eo[\hat{f}]-\Eo[f_{\H}]\le r(l,N,\lambda)$ bounds, 
	where $\lambda$ is a regularization parameter converging to zero as we see more samples ($l\rightarrow \infty$, $N=N_i\rightarrow \infty$), and
	choose the $(l,N,\lambda)$ triplet appropriately to drive $r(l,N,\lambda)$ and hence $\Eo[\hat{f}]-\Eo[f_{\H}]$ to $0$.
  \item obtaining convergence rates: We establish convergence rates for a general prior family $\P(b,c)$ \cite{caponnetto07optimal}, where $b$ captures the effective input dimension, and larger $c$ means smoother $f_{\H}$. 
		  In particular, when  $l = N^a$ ($a > 0$), the effective dimension is small (large $b$), and the total number of samples processed $t=lN=N^{a+1}$ is fixed, 
		  one obtains a rate of $1/t^{2/7}$ for a smooth regression function ($c=2$), $1/t^{1/5}$ in the non-smooth case ($c=1$). %In other words, for a given total computation budget $t$, we prove that convergence occurs at a faster rate for smoother problems.
\end{compactenum}
The motivation for considering the  $\P(b,c)$ family is two-fold:
		  \begin{compactenum}
			\item it does not assume parametric distributions, still certain complexity terms can be explicitly upper bounded in the family. This property will be exploited in our analysis.
			\item (for special input distributions) parameter $b$ can be related to the spectral decay of Gaussian Gram matrices, thus available analysis techniques \cite{steinwart08support} might give alternative prior characterizations.
		  \end{compactenum}

Briefly, we focus on the following question:
\begin{center}
  \fbox{\begin{minipage}{18.5em}
  Can the distribution regression problem be solved consistently under mild conditions?
  \end{minipage}}
\end{center}
Despite the large number of available ``solutions'' and applications of distribution regression dating back to $1999$ \cite{haussler99convolution}, surprisingly this pretty fundamental question has hardly
 been touched. In our paper we give affirmative answer to the question by presenting the analysis of a \emph{simple} kernel ridge regression approach 
[see Eq.~\eqref{eq:MERR}] in the two-stage sampled 
($\M \rightarrow \b{z} \rightarrow \hat{\b{z}}$) setting. 

{\bf Review of approaches to learning on distributions:} A number of methods have been proposed over the years to compute the similarity of distributions or bags of samples. As a first approach, one could fit
a parametric model to the bags, and estimate the similarity of the bags based on the obtained parameters. It is then possible to
define learning algorithms on the basis of these similarities, which often take analytical form. Typical examples with explicit formulas include Gaussians,
finite mixtures of Gaussians, and distributions from the exponential family (with known log-normalizer function and zero carrier measure) \cite{wang09closed,kondor03kernel,jebara04probability,nielsen12closed}.
A major limitation of these methods, however, is that they apply quite simple parametric assumptions, which may not be sufficient or verifiable in practise.

A heuristic related to the parametric approach is to assume that the training distributions are Gaussians in a reproducing kernel Hilbert space; see for
example \cite{jebara04probability,zhou06sample} and references therein. This assumption is algorithmically appealing, as many divergence measures for Gaussians can be
computed in closed form using only inner products, making them straightforward to kernelize. A fundamental shortfall of kernelized Gaussian divergences is the lack of
their consistency analysis in specific learning algorithms.

A more theoretically grounded approach to learning on distributions has been to define positive definite kernels \cite{scholkopf02learning} on the basis of statistical divergence measures on distributions,
or by metrics on non-negative numbers; these can then be used in kernel algorithms.
This category includes work on semigroup kernels \cite{cuturi05semigroup}, nonextensive information theoretical kernel constructions \cite{martins09nonextensive}, and kernels based
on Hilbertian metrics \cite{hein05hilbertian}. For example, in \cite{cuturi05semigroup} the intuition is as follows: if two measures or sets of points
overlap, then their sum is expected to be more concentrated. The value of dispersion can be measured by entropy or inverse generalized variance. In the second type of approach
\cite{hein05hilbertian}, homogeneous Hilbert metrics on the non-negative real line are used to define the similarity of probability distributions.
While these techniques guarantee to provide valid kernels on certain restricted domains of measures, the performance of learning
algorithms based on finite sample estimates of these kernels remains a challenging open question.
One might also plug into learning algorithms (based on similarities of distributions) consistent R{\'e}nyi and Tsallis divergence estimates \cite{poczos11nonparametric,poczos12support}, but
these similarity indices are \emph{not} kernels, and their consistency in specific learning tasks, similarly to the previous works, is open.

To the best of our knowledge, the only prior work addressing the consistency of regression on distributions requires kernel density estimation
\cite{poczos13distribution,oliva14fast}, assumes that the response variable is scalar-valued\footnote{\cite{oliva13ICML} considers the case where the responses are also distributions.}, and the covariates are nonparametric continuous distributions on $\R^d$.
As in our setting, the exact forms of these distributions are unknown; they are available only through finite sample sets. P{\'o}czos et al.\ estimated these distributions
through a kernel density estimator (assuming these distributions to have a density) and then constructed a kernel regressor that acts on these kernel
density estimates.\footnote{We would like to clarify that the kernels used in their work are classical smoothing kernels (extensively studied in non-parametric statistics \cite{gyorfi2002})
and not the reproducing kernels that appear throughout our paper.} Using the classical bias-variance decomposition analysis for kernel regressors, they show the consistency of the 
constructed kernel regressor, and provide a polynomial upper bound on the rates, assuming the true regressor to be H{\"o}lder continuous, and the meta distribution that generates the 
covariates $x_i$ to have finite doubling dimension \cite{kpotufe2011}.\footnote{Using a random kitchen sinks approach, with orthonormal basis 
projection estimators and RBF kernels \cite{oliva14fast} proposes a distribution regression algorithm that can computationally handle large scale datasets; 
as with \cite{poczos13distribution}, this approach is based on density estimation in $\R^d$.\label{footnote:RNDkitchensink}}

An alternative paradigm in learning when the inputs are ``bags of objects'' is to simply treat each input as a finite
set, and to define kernel learning algorithms based on set kernels \cite{gartner02multi} (also called multi-instance kernels or ensemble kernels, and instances of convolution kernels \cite{haussler99convolution}). In this case, the similarity of two sets is measured by the average pairwise point
similarities between the sets.
From a \emph{theoretical} perspective,
very little has been done to establish the consistency of set kernels in learning since their introduction in 1999 \cite{haussler99convolution,gartner02multi}:
i.e.\ in what sense (and with what rates) is the learning algorithm consistent,
 when the number of items per bag, and the number of bags, is allowed to increase?

It is possible, however, to view set kernels in a distribution setting, as they represent valid kernels between (mean) embeddings of empirical probability measures into a reproducing kernel Hilbert space (RKHS) \cite{BerTho04}.
The \emph{population limits are well-defined} as being dot products between the embeddings of the generating distributions \cite{altun06unifying}, and for characteristic kernels the distance between embeddings defines a \emph{metric} on probability measures \cite{gretton12kernel,sriperumbudur11universality}.  When bounded kernels are used, mean embeddings \emph{exist for all probability measures} \cite{fukumizu04dimensionality}.
When we consider the distribution regression setting, however, there is no reason to limit ourselves to set kernels.
Embeddings of probability measures to RKHS are  used by \cite{christmann10universal} in defining a yet larger class
of easily computable kernels on distributions, via operations performed on the embeddings and their distances.
Note that the relation between set kernels and kernels on distributions has been applied by \cite{muandet12learning} for classification on distribution-valued inputs, however consistency was not studied in that work.

Our \tb{contribution} in this paper is to establish the consistency of an algorithmically simple, mean embedding based ridge regression method
(described in Section \ref{sec:problem}) for the distribution regression problem.
This result applies both to the basic set kernels of \cite{haussler99convolution,gartner02multi}, the distribution kernels of  \cite{christmann10universal}, and additional related kernels proposed herein.
We provide two-stage sampled excess error bounds, consistency proof and convergence rates in Section \ref{sec:convergence analysis}, and break down the various tradeoffs arising in
different sample size and problem difficulties. The principal challenge in proving theoretical guarantees arises from the two-stage sampled nature of the inputs. 
In our analysis, we make use of \cite{caponnetto07optimal}, who 
provide error bounds for the one-stage sample setup. These results will make our analysis somewhat shorter (but still rather challenging) by giving upper bounds for some of the upcoming objective terms. 
Even the verification of these conditions requires care (Section \ref{sec:assumptions}) since the inputs in the ridge regression are themselves distribution embeddings (i.e., functions in a reproducing kernel Hilbert space).

Due to the differences in the assumptions made and the loss function used, a direct comparison of our theoretical result and that of 
\cite{poczos13distribution}\textsuperscript{\ref{footnote:RNDkitchensink}} remains an open question, 
however we make two observations. First, our approach is more general, since we may regress from any probability measure defined on a separable, topological domain endowed with a kernel.
P{\'o}czos et al.'s work is restricted to compact domains of finite dimensional Euclidean spaces, and requires the distributions to admit probability 
densities; distributions on strings, time series, 
graphs, and other structured objects are disallowed. Second, density estimates in high dimensional spaces suffer from slow convergence rates \cite[Section 6.5]{Wasserman06}. Our approach 
avoids this problem, as it works directly on distribution embeddings, and does not make use of density estimation as an intermediate step.

\vspace*{-0.1cm}
\section{THE DISTRIBUTION REGRESSION PROBLEM} \label{sec:problem}
\vspace*{-0.1cm}
In this section, we define the distribution regression problem, for a general RKHS on distributions. In Section~\ref{sec:assumptions}, we will provide examples of valid kernels for this RKHS,  including  set kernels  \cite{haussler99convolution,gartner02multi}, the kernels from \cite{christmann10universal}, and further related kernels. Below, we first introduce some notation and then formally discuss the distribution regression problem. 

\textbf{Notation:} Let $(\X,\tau)$ be a topological space and let $\Bo(\X):=\Bo(\tau)$ be the Borel $\sigma$-algebra induced by the topology $\tau$. $\M_1^+(\X)$ denotes the set of Borel probability measures on $(\X,\Bo(\X))$. 
The weak topology ($\tau_w=\tau_w(\X,\tau)$) on $\M^+_1(\X)$ is defined as the weakest topology such that the $L_h:(\M^+_1(\X),\tau_w) \rightarrow \R$, 
$L_h(x)=\int_{\X}h(u)\d x(u)$ mapping is continuous for all $h\in C_b(\X)=\{(\X,\tau)\rightarrow  \R\text{ bounded, continuous functions}\}$.
Let $H=H(k)$ be the RKHS \cite{steinwart08support} with $k:\X\times \X \rightarrow \R$ as the reproducing kernel. 
Denote by 
\begin{align*}
  X &= \mu\left(\M^+_1\left(\X\right)\right) = \{\mu_x: x\in \M^+_1\left(\X\right)\}\subseteq H 
 \end{align*}
the set of $\mu_x =\int_{\X}k(\cdot,u)\d x(u)=\E_{u\sim x}[k(\cdot,u)] \in H$ mean embeddings \cite{BerTho04} of the distributions 
to the space $H$, and let $Y=\R$. Intuitively, $\mu_x$ is the canonical feature map [$k(\cdot,u)$] averaged according to the probability measure
[$\d x(u)$]. Let $\H=\H(K)$ be the RKHS of functions with $K:X\times X\rightarrow \R$ as the reproducing kernel. 
$\L(\H)$ is the space of $\H\rightarrow\H$ bounded linear operators, and $\delta_{\mu_a}$ denotes the evaluation functional at $\mu_a$ ($a\in\M^+_1(\X)$). 
For $M\in\L(\H)$ the operator norm is defined as $\left\|M\right\|_{\L(\H)}=\sup_{0\ne q\in \H}\frac{\left\|Mq\right\|_{\H}}{\left\|q\right\|_{\H}}$.
Given $(U_1,\S_1)$ and $(U_2,\S_2)$ measurable spaces the $\S_1\otimes \S_2$ product $\sigma$-algebra \cite[page~480]{steinwart08support} on the product space 
$U_1 \times U_2$ is the $\sigma$-algebra generated by the cylinder sets $U_1\times S_2$, $S_1\times U_2$ ($S_1\in \S_1$, $S_2\in \S_2$). $\E[\cdot]$ denotes expectation.

\textbf{Distribution regression:} In the  \emph{distribution regression} problem, we are given samples $\hat{\b{z}} = \{(\{x_{i,n}\}_{n=1}^{N_i},y_i)\}_{i=1}^l$ with $x_{i,1},\ldots, x_{i,N_i} \stackrel{i.i.d.} {\sim} x_i$ where $\b{z}=\{(x_i,y_i)\}_{i=1}^l$ with $x_i\in \M^+_1\left(\X\right)$ and $y_i\in Y$  drawn i.i.d.~from a joint meta distribution 
$\mathcal{M}$ defined on the measurable space $(\M^+_1(\X)\times\R,\Bo(\tau_w)\otimes \Bo(\R))$. Unlike in classical supervised learning problems, the problem at hand involves two levels of randomness, wherein first $\b{z}$ is drawn from $\mathcal{M}$ and then $\hat{\b{z}}$ is generated by sampling points from $x_i$ for all $i=1,\ldots,l$. The \tb{goal} is to learn the relation between the random distribution $x$ and scalar response $y$ based on 
the observed $\hat{\b{z}}$. For notational simplicity, we will assume that $N=N_i$ ($\forall i$).

As in the classical regression task ($\R^d\rightarrow \R$), distribution regression can be tackled as a kernel ridge regression problem (using squared loss as the discrepancy criterion). The 
kernel (say $\Kmeta$) is defined on $\M^+_1(\X)$, and the regressor is then modelled  by an element in the RKHS $\Hmeta=\Hmeta(\Kmeta)$ of functions mapping from $\M^+_1(\X)$ to $\R$. In this paper, we choose 
$\Kmeta(x,x')=K(\mu_x,\mu_{x'})$ where $x,x'\in \M^+_1(\X)$ and so that 
the function (in $\Hmeta$) to describe the $(x,y)$ random relation is constructed as a composition
\begin{align*}
  \M^+_1\left(\X\right) \xrightarrow{\mu} X (\subseteq H=H(k)) \xrightarrow{f\in \H=\H(K)} \R.
\end{align*}
In other words, the distribution $x\in \M^+_1\left(\X\right)$ is first mapped to $X\subseteq H$ by the mean embedding $\mu$, and the result is composed with $f$, 
an element of the RKHS $\H=\H(K)$. Assuming that $f_{\Hmeta}$, the minimizer of the expected risk ($\Eo$) over $\Hmeta$ exists, then a function  $f_{\H}$ also exists, and satisfies
\begin{align*}
    \Eo\left[f_{\H}\right]&\hspace{-0.1cm} = \hspace{-0.1cm}\inf_{f\in \H} \Eo[f]\hspace{-0.1cm} = \hspace{-0.1cm}\inf_{f\in \H} \E_{(x,y)\sim\mathcal{M}} [f(\mu_x) - y]^2\\
			  &\hspace{-0.1cm} = \hspace{-0.1cm}\inf_{g\in \Hmeta}\E_{(x,y)\sim\mathcal{M}}[g(x)-y]^2\hspace{-0.1cm} =\hspace{-0.1cm} \inf_{g\in \Hmeta} \Eo[g] = \Eo\left[g_{\Hmeta}\right].
\end{align*}
The classical regularization approach is to optimize
\begin{align}
	f_{\b{z}}^{\lambda} &= \argmin_{f\in \H}\frac{1}{l}\sum_{i=1}^l\left[f(\mu_{x_i})-y_i\right]^2 + \lambda \left\|f\right\|_{\H}^2\label{eq:obj0}
\end{align}
 instead of $\Eo$, based on samples $\b{z}$. Since $\b{z}$ is not accessible, we consider the objective function defined by 
the observable quantity $\hat{\b{z}}$,
\begin{align}
    f_{\hat{\b{z}}}^{\lambda} &= \argmin_{f\in \H}\frac{1}{l}\sum_{i=1}^l\left[f(\mu_{\hat{x}_i}) - y_i\right]^2 + \lambda \left\|f\right\|_{\H}^2, \label{eq:obj}
\end{align}
where $\hat{x}_i=\frac{1}{N}\sum_{n=1}^N \delta_{x_{i,n}}$ is the empirical distribution determined by $\left\{x_{i,n}\right\}_{i=1}^N$. Algorithmically, ridge regression is quite 
simple \cite{cucker02mathematical}: given training samples $\hat{\b{z}}$, the prediction for a new $t$ test distribution is
\begin{empheq}[box=\fbox]{align}
    (f_{\hat{\b{z}}}^{\lambda} \circ \mu) (t) & = [y_1,\ldots,y_l](\b{K}+l \lambda \b{I}_l)^{-1} 
   \b{k}\in \R,\label{eq:MERR}\\
  \b{K} &=  [K(\mu_{\hat{x}_i},\mu_{\hat{x}_j})]\in \R^{l\times l},\nonumber\\
  \b{k} &= \left[K(\mu_{\hat{x}_1},\mu_t); \ldots; K(\mu_{\hat{x}_l},\mu_t)\right]\in \R^l.\nonumber
\end{empheq}

\tb{Remarks:} 
\begin{compactenum}
    \item It is important to note that the algorithm has access to the sample points \emph{only via} their \emph{mean embeddings} $\{\mu_{\hat{x}_i}\}_{i=1}^l$ in Eq.~\eqref{eq:obj}.
    \item There is a \emph{two-stage sampling difficulty} to tackle: The transition from $f_{\H}$ to $f_{\b{z}}^{\lambda}$ represents the fact that we have only $l$ distribution samples ($\b{z}$); 
	the transition from $f_{\b{z}}^{\lambda}$ to $f_{\hat{\b{z}}}^{\lambda}$ means that the $x_i$ distributions can be accessed only via samples 
	($\hat{\b{z}}$).
    \item While  ridge regression can be performed using the kernel $\Kmeta$, the two-stage sampling makes it difficult to work with arbitrary $\Kmeta$. By contrast, our 
	  choice of $\Kmeta(x,x')=K(\mu_x,\mu_{x'})$ enables us to handle the two-stage sampling by estimating $\mu_x$ with an empirical estimator and using it in the algorithm as shown above.
\end{compactenum}

The main \tb{goal} of this paper is to analyse the excess risk
$\Eo[f_{\hat{\b{z}}}^{\lambda}]-\Eo[f_{\H}]$, i.e., the regression performance compared to the best possible estimation from $\H$, and to establish consistency and rates of convergence as a function of the
$(l,N,\lambda)$ triplet, and of the difficulty of the problem in the sense of \cite{caponnetto07optimal}.

\vspace*{-0.1cm}
\section{ASSUMPTIONS} \label{sec:assumptions}
\vspace*{-0.1cm}
In this section we detail our assumptions on the $(\X,k,K)$ triplet, and show that regressing with set kernels fit into the studied problem family. 
Our analysis will rely on  existing ridge regression results \cite{caponnetto07optimal} 
which focus on problem \eqref{eq:obj0}, where only a single-stage sampling is present; hence we have to verify the associated conditions.
Though we make use of these results, the analysis still remains rather challenging; the available bounds can moderately shorten our proof.
We must also take particular care in verifying that \cite{caponnetto07optimal}'s conditions are met, since they must hold for the space of \emph{mean embeddings of the distributions} ($X=\mu\left(\M_1^+(\X)\right)$), whose properties as a function of $\X$ and $H$ must themselves be established.
Our assumptions:
\begin{compactitem}
  \item $\exists f_{\H}$ such that $\Eo[f_{\H}]=\inf_{f\in \H} \Eo(f)$.
  \item $(\X,\tau)$ is a separable, topological domain. 
  \item $k: \X \times \X\rightarrow \R$ is bounded ($\exists B_k<\infty$ such that $\sup_{u\in\X}k(u,u)\le B_k$) and continuous.
  \item $K: X \times X \rightarrow \R$ is bounded, i.e., $\exists B_K<\infty$ such that
		      \begin{align}
			  K(\mu_a,\mu_a) &\le B_K, \quad (\forall \mu_a\in X), \label{eq:bounded kernel}
		      \end{align} and $\Psi(\mu_c) := K(\cdot,\mu_c):X\rightarrow \H$ is H{\"o}lder continuous, i.e., $\exists L>0$, $h\in (0,1]$ such that for $\forall (\mu_a,\mu_b)\in X\times X$
      	        \begin{align}
		      \left\|\Psi(\mu_a) - \Psi(\mu_b)\right\|_{\H} &\le L \left\|\mu_a - \mu_b\right\|_H^h.\label{eq:K:Lip}
		 \end{align}
  \item $y$ is bounded: $\exists C<\infty$ such that $|y|\le C$ almost surely.
  \item $X=\mu(\M^+_1(\X)) \in \Bo(H)$.
\end{compactitem}

\tb{Discussion of the assumptions:} We give a short insight into the consequences of our assumptions and present some concrete examples.
\begin{compactenum}
    \item The boundedness and continuity of $k$ imply the measurability of $\mu: (\M^+_1(\X),\Bo(\tau_w))\rightarrow (H,\Bo(H))$, which using the $X\in\Bo(H)$ condition guarantees that 
	  the $\rho$, the measure induced by $\M$ on $X\times \R$ is well-defined (see the supplementary material).
    \item For a linear kernel, $K(\mu_a,\mu_b)=\left<\mu_a,\mu_b\right>_H$, $(\mu_a,\mu_b\in X)$, one can verify (see the supplementary material) that H{\"o}lder continuity holds with $L=1$, $h=1$. 
      Also, since $K(\mu_a,\mu_b)\le B_k$ for any $a,b\in \M^+_1(\X)$, we can choose $B_K=B_k$. Evaluating the kernel, $K$ at the $\mu_{\hat{x}_i} = \int_{\X}k(\cdot,u)\d \hat{x}_i(u) =  \frac{1}{N}\sum_{n=1}^N k(\cdot,x_{i,n})$ empirical embeddings yields the standard set kernel:% by the bilinearity of $\left<\cdot,\cdot\right>_H$ and the reproducing property of $k$:
      \begin{align*}
	    K(\mu_{\hat{x}_i},\mu_{\hat{x}_j}) &= \frac{1}{N^2}\sum_{n,m=1}^N k(x_{i,n},x_{j,m}).
	\end{align*}
  \item One can also prove (see the supplement) by using the properties of negative/positive definite functions \cite{berg84harmonic} 
		  that many $K$ functions on $X\times X$ are kernels and (in case of compact metric $\X$ domains) H{\"o}lder continuous.\footnote{To guarantee the H{\"o}lder property of $K$-s, we assume the continuity of $\mu$.
		  For example, if $\X$ is a compact metric space and $k$ is universal, then $\mu$ metrizes the weak topology $\tau_w$ \cite[Theorem~23, page~1552]{sriperumbudur10hilbert}, hence $\mu$ is continuous. In this case $X=\mu(\M^+_1(\X))$ 
		  is compact metric (see the supplement), 
		  thus closed and hence $X\in \Bo(H)$ also holds.} Some examples are listed in Table~\ref{tab:K examples}; these kernels are the natural 
		  extensions to distributions of the Gaussian \cite{christmann10universal}, exponential, Cauchy, generalized t-student and inverse multiquadratic kernels.
    \item $Y=\R$ is a separable Hilbert space hence Polish, and thus the $\rho(y|\mu_a)$ conditional distribution ($y\in\R$, $\mu_a\in X$) is well-defined; see  \cite[Lemma~A.3.16, page~487]{steinwart08support}. 
    \item The separability of $\X$ and the continuity of $k$ implies the separability of $H$ \cite[Lemma 4.33, page~130]{steinwart08support}. Also, since $X\subseteq H$, $X$ is separable; hence so is $\H$ due to the continuity of $K$.
\end{compactenum}

\tb{Verification of \cite{caponnetto07optimal}'s conditions:} Below we prove that \cite{caponnetto07optimal}'s conditions hold under our assumptions.
\begin{compactenum}
  \item $Y=\R$ and $\H$ are separable Hilbert spaces -- as we have seen.
  \item By the bilinearity of $\left<\cdot,\cdot\right>_{\H}$ and the reproducing property of $K$, the measurability of 
	    $(\mu_x,\mu_t)\mapsto \left<K(\cdot,\mu_x)w,K(\cdot,\mu_t)v\right>_{\H} 
			 =w K(\mu_x,\mu_t) v \quad (\forall w,v\in \R)$
	is equivalent to that of $(\mu_x,\mu_t)\mapsto K(\mu_x,\mu_t)$; the latter follows from the H{\"o}lder continuity of $\Psi$ (see the supplement).
  \item Due to the boundedness of $y$, we have
	  $\int_{X\times \R} y^2 \d \rho(\mu_x,y)  \le     \int_{X\times \R} C^2 \d \rho(\mu_x,y) \le C^2<\infty$,
    and $\exists \Sigma>0, \exists M>0$ such that 
    \begin{align}
	  \hspace*{-0.25cm}\int_{\R}e^{\frac{\left|y-f_{\H}(\mu_x)\right|}{M}}-\frac{\left|y-f_{\H}(\mu_x)\right|}{M}-1 \d  \rho(y|\mu_x) &\le
	  \frac{\Sigma^2}{2M^2} \label{eq:bounded noise}
    \end{align}
    for $\rho_X\text{-almost }\mu_x\in X$, where $\rho(\mu_x,y) = \rho(y|\mu_x)\rho_X(\mu_x)$ is factorized into conditional and marginal distributions.
    \eqref{eq:bounded noise} is a model of the noise of the output $y$; it is satisfied, for example in case of bounded noise \cite[page~9]{caponnetto07optimal}.
    By the boundedness of $y$ and that of kernel $K$ this property holds:
    $|y-f_{\H}(\mu_x) | \le |y| +  |f_{\H}(\mu_x)| \le  C + \left\|f_{\H}\right\|_H \sqrt{B_K}$,
    where we used the triangle inequality and Lemma~4.23 (page~124) from \cite{steinwart08support}.
\end{compactenum}
  
\begin{table*}
  %\small
  \begin{center}
  \caption{Nonlinear kernels on mean embedded distributions: $K=K(\mu_a,\mu_b)$; $\theta>0$. For the H{\"o}lder continuity, we assume that $\X$ is a compact metric space and $\mu$ is continuous (the latter is implied e.g., by a universal $k$).}
  \label{tab:K examples}
  \begin{tabular}{@{}ccccc@{}}
    \toprule
      $K_G$ & $K_e$ & $K_C$ & $K_t$ & $K_i$\\\midrule
      $e^{-\frac{\left\|\mu_a-\mu_b\right\|_H^2}{2\theta^2}}$ & $e^{-\frac{\left\|\mu_a-\mu_b\right\|_H}{2\theta^2}}$ & $\left(1+\left\|\mu_a-\mu_b\right\|_H^2 / \theta^2\right)^{-1}$ & $\left(1+\left\|\mu_a-\mu_b\right\|_H^{\theta}\right)^{-1}$ & $\left(\left\|\mu_a-\mu_b\right\|_H^2+\theta^2\right)^{-\frac{1}{2}}$\\
      $h=1$ & $h=\frac{1}{2}$ & $h=1$  & $h=\frac{\theta}{2}$ ($\theta\le 2$) & $h=1$\\ \bottomrule
  \end{tabular}
  \end{center}
\end{table*}

\vspace*{-0.1cm}
\section{ERROR BOUNDS, CONSISTENCY, CONVERGENCE RATE} \label{sec:convergence analysis}
\vspace*{-0.1cm}
In this section, we present our main result: we
derive high probability upper bound for the excess risk $\Eo\left[f^{\lambda}_{\hat{\b{z}}}\right] - \Eo\left[f_{\H}\right]$ of the 
mean embedding based ridge regression (MERR) method, see our main 
theorem. We also illustrate the upper bound for particular classes of prior distributions, resulting in sufficient conditions for convergence and concrete convergence rates (see Consequences~\ref{conseq:excess-rate}-\ref{conseq:conv-rate}).
We first give a high-level sketch of our convergence analysis and the results are stated with their intuitive interpretation. Then an outline of the main proof ideas follows; technical details of the proof steps may be found in the 
supplement.

At a high level, our convergence analysis takes the following form: Having explicit expressions for $f_{\b{z}}^{\lambda}$, $f_{\hat{\b{z}}}^{\lambda}$ 
[see Eq.~\eqref{eq:f_zlambda}-\eqref{eq:f_hatz^lambda}], we will decompose
the excess risk $\Eo[f^{\lambda}_{\hat{\b{z}}}] - \Eo[f_{\H}]$ into five terms:
{
\setlength{\abovedisplayskip}{11pt}
\setlength{\belowdisplayskip}{11pt}
\setlength{\abovedisplayshortskip}{11pt}
\setlength{\belowdisplayshortskip}{11pt}
\begin{align*}
		\Eo\left[f^{\lambda}_{\hat{\b{z}}}\right] - \Eo\left[f_{\H}\right] &\le 5 \left[S_{-1} + S_0 + \A(\lambda) + S_1 + S_2 \right],\\
	      S_{-1} &= \|\sqrt{T} (T_{\hat{\b{x}}}+\lambda)^{-1}(g_{\hat{\b{z}}} - g_{\b{z}})\|_{\H}^2, \\
	      S_0 &= \|\sqrt{T} (T_{\hat{\b{x}}}+\lambda)^{-1}(T_{\b{x}}-T_{\hat{\b{x}}}) f^{\lambda}_{\b{z}}\|_{\H}^2,\\
	      S_1 &= \|\sqrt{T}(T_{\b{x}}+\lambda)^{-1}(g_{\b{z}}-T_{\b{x}}f_{\H})\|_{\H}^2,\\
	      S_2&= \|\sqrt{T}(T_{\b{x}}+\lambda)^{-1}(T-T_{\b{x}})(f^{\lambda}-f_{\H})\|_{\H}^2,\\
	      \A(\lambda) &= \|\sqrt{T}(f^{\lambda}-f_{\H})\|_{\H}^2, \vspace*{-1cm}
\end{align*}
}
where $f^{\lambda} = \argmin_{f\in\H}\Eo[f]+\lambda\left\|f\right\|_{\H}^2$, $T_{\mu_a} = K(\cdot,{\mu_a})\delta_{\mu_a}$ [$T_{\mu_a}(f)=K(\cdot,\mu_a)f(\mu_a)$, $\mu_a\in X$],
{
\setlength{\abovedisplayskip}{11pt}
\setlength{\belowdisplayskip}{11pt}
\setlength{\abovedisplayshortskip}{11pt}
\setlength{\belowdisplayshortskip}{11pt}
\begin{align}
  T = \int_{X}T_{\mu_a}\d \rho_X(\mu_a)\in \L(\H), \quad T_{\mu_a}\in \L(\H). \label{eq:Tdef}
\end{align}
}
\begin{compactenum}
    \item Three of the terms ($S_1$, $S_2$, $\A(\lambda)$) will be identical to terms in \cite{caponnetto07optimal}, hence their bounds can be applied.
    \item The two new terms ($S_{-1}$, $S_0$), the result of two-stage sampling, will be upper bounded by making use of the convergence of the empirical mean embeddings, and the H{\"o}lder property of $K$. 
\end{compactenum}
These bounds will lead to the following results:
	      \begin{theorem}\label{theo1}
		    Let $M$, $\Sigma$  and $T$ be as in \eqref{eq:bounded noise}, \eqref{eq:Tdef}.
		    Let $\Psi(\mu_a) = K(\cdot,\mu_a): X\rightarrow \H$ be H{\"o}lder continuous with constants $L$, $h$.
		    Let $l\in\mathbb{N}$, $N\in\mathbb{N}$, $\lambda>0$,  $0<\eta<1$, $C>0$, $\delta>0$, $C_{\eta}=32\log^2(6/\eta)$, $|y|\le C$ (a.s.), $\A(\lambda)$ the residual as above, and define
		    $\B(\lambda)=\|f^{\lambda} - f_{\H}\|_{\H}^2$ the reconstruction error, $\N(\lambda) = Tr [(T+\lambda)^{-1}T]$ the effective dimension.
     		Then with probability at least $1-\eta-e^{-\delta}$
		    \begin{eqnarray*}
			\lefteqn{\Eo\left[f^{\lambda}_{\hat{\b{z}}}\right] - \Eo\left[f_{\H}\right]\le}&\\
			&&\hspace*{-0.7cm}\le 5\Bigg\{ \frac{4L^2 C^2 \left(1+\sqrt{\log(l)+\delta}\right)^{2h} (2B_k)^h}{\lambda N^h} \left[  1 + \frac{4(B_K)^2}{\lambda^2}\right] + \\
			&&\hspace*{-0.6cm}   \A(\lambda)+C_{\eta} \left[\frac{B_K^2\B(\lambda)}{l^2\lambda} + \frac{B_K \A(\lambda)}{4l\lambda} + \frac{B_K M^2}{l^2\lambda} + \frac{\Sigma^2 \N(\lambda)}{l} \right] \Bigg\}
		    \end{eqnarray*}
		    provided that $l\ge 2C_{\eta}B_K\N(\lambda) / \lambda$, $\lambda \le \left\|T\right\|_{\L(\H)}$, $N \ge (1+\sqrt{\log(l)+\delta})^{2} 2^{\frac{h+6}{h}}B_k (B_K)^{\frac{1}{h}} L^{\frac{2}{h}} / \lambda^{\frac{2}{h}}$.
	      \end{theorem}
	      Below we specialize our bound on the excess risk  for a general prior class, which captures the difficulty of the regression problem as defined in \cite{caponnetto07optimal}. This 
    	      $\P(b,c)$ class is described by two parameters $b$ and $c$: intuitively, larger $b$ means faster decay of the eigenvalues of the covariance operator $T$ 
	      [\eqref{eq:Tdef}], hence smaller effective input dimension; larger $c$ corresponds to smoother $f_{\H}$. Formally:
	
	      \tb{Definition of the $\P(b,c)$ class:} Let us fix the positive constants $M$, $\Sigma$, $R$, $\alpha$, $\beta$. Then given $1<b$, $c\in[1,2]$, the $\P(b,c)$ class is
	      the set of probability distributions $\rho$ on $Z=X\times \R$  such that (i)
              the $(\mu_x,y)$ assumption holds with $M$, $\Sigma$ in \eqref{eq:bounded noise}, 
	      (ii) there is a $g\in\H$ such that $f_{\H}=T^{\frac{c-1}{2}}g$ with $\left\|g\right\|_{\H}^2\le R$, (iii)
	      in the $T=\sum_{n=1}^Nt_n\left<\cdot,e_n\right>_{\H}e_n$ spectral theorem based decomposition ($(e_n)_{n=1}^N$ is a basis of $ker(T)^{\perp}$), $N=+\infty$, and the eigenvalues of 
	      $T$ satisfy $\alpha \le n^b t_n\le \beta \quad (\forall n\ge 1)$.

 	      We can provide a simple example of when the source decay conditions hold, 
   	      in the event that the distributions are normal with means $m_i$ and identical variance ($x_i = N(m_i,\sigma^2I$)).
	      When Gaussian kernels ($k$) are used with linear $K$, 
	      then $K(\mu_{x_i},\mu_{x_j}) = e^{-c\left\|m_i-m_j\right\|^2}$ \cite[Table 1, line 2]{muandet12learning}  
	      (Gaussian, with arguments equal to the difference in means). Thus, this Gram matrix will correspond to the Gram matrix using a Gaussian kernel between points 
 	      $m_i$. The spectral decay of the Gram matrix will correspond to that of the Gaussian kernel, with points drawn from the meta-distribution over the $m_i$. 
	      Thus, the source conditions are analysed in the same manner as for Gaussian Gram matrices, e.g. see \cite{steinwart08support} for a discussion of the spectral decay properties.

	      In the $\P(b,c)$ family, the behaviour of $\A(\lambda)$, $\B(\lambda)$ and $\N(\lambda)$ is known; specializing our theorem we get:\footnote{In what follows, we assume the conditions of the main theorem and $\rho\in\P(b,c)$.}
	      \begin{consequence}[Excess risk in the $\P(b,c)$ class]\label{conseq:excess-rate} 
	      \begin{eqnarray*}
		  \hspace*{0.2cm}\lefteqn{\Eo\left[f^{\lambda}_{\hat{\b{z}}}\right] - \Eo\left[f_{\H}\right]\le 5 \Bigg\{\frac{4L^2 C^2 \left(1+\sqrt{\log(l)+\delta}\right)^{2h} (2B_k)^h}{\lambda N^h}} \\
		  &&\hspace*{-0.7cm}   \times\left[  1 + \frac{4(B_K)^2}{\lambda^2}\right] + R\lambda^c + C_{\eta}\times\\
		  &&\hspace*{-0.7cm}   \left[\frac{B_K^2R\lambda^{c-2}}{l^2}+\frac{B_K R\lambda^{c-1}}{4l}+\frac{B_K M^2}{l^2\lambda}+\frac{\Sigma^2\beta b}{(b-1)l\lambda^{\frac{1}{b}}}\right]\Bigg\}.
	      \end{eqnarray*}	      
	      \end{consequence}
	      By choosing $\lambda$ appropriately as a function of $l$ and $N$, the excess risk $\Eo[f^{\lambda}_{\hat{\b{z}}}] - \Eo[f_{\H}]$ converges to $0$, and we can use
	      Consequence~1 to  obtain convergence rates: the task reduces to the study of
	      \begin{align}
		  r(l,N,\lambda) &= \frac{\log^h(l)}{N^h\lambda^3}+\lambda^c + \frac{1}{l^2\lambda} + \frac{1}{l\lambda^{\frac{1}{b}}}\rightarrow 0,\label{eq:r}
	      \end{align}
	      subject to $l\ge \lambda^{-\frac{1}{b}-1}$.\footnote{Note that the $N\ge \log(l) / \lambda^{\frac{2}{h}}$ constraint has been discarded; it is implied by the convergence of the first term in $r$ [Eq.~\eqref{eq:r}] (see the supplementary material).} By matching two terms in \eqref{eq:r}, solving  for $\lambda$ and plugging the result back to the bound (see the supplementary material), we obtain:
	      \begin{consequence}[Consistency and convergence rate in $\P(b,c)$]\label{conseq:conv-rate} Let $l=N^a$ ($a>0$). The
		excess risk can be upper bounded (constant multipliers are discarded) by the quantities given in the last column of Table~\ref{tab:conv-rates}.
	      \end{consequence}
	      \tb{Note:} in function $r$ [Eq.~\eqref{eq:r}] (i) the first term comes from the error of the mean embedding estimation, (ii) the second term corresponds to $\A(\lambda)$, a
	      complexity measure of $f_{\H}$, (iii) the third term is due to the $S_1$ bound, (iv) the fourth term expresses $\N(\lambda)$, a complexity index of the hypothesis space $\H$ according to the marginal measure $\rho_X$.
	      As an example, let us take two rows from Table~\ref{tab:conv-rates}:
	      \begin{compactenum}
		  \item
			First row: In this case the first and second terms dominate $r(l,N,\lambda)$ in \eqref{eq:r}; in other words the error is determined by the mean embedding estimation
			process and the complexity of $f_{\H}$. Let us assume that $b$ is large in the sense that $1/b\approx 0$, $(b+1)/b\approx 1$ (hence, the effective dimension of the input space is small); %AG
			and assume that $K$ is Lipschitz ($h=1$). Under these conditions the lower bound for $a$ is approximately $\max (c/(c+3),1/(c+3)) = c/(c+3)\le a$
			(since $c\ge 1$). Using such an $a$ (i.e., the exponent in $l=N^a$ is not too small), then the convergence rate is $\left[\log(N)/N\right]^{\frac{c}{c+3}}$. Thus, for example, if $c=2$ ($f_{\H}=T^{\frac{c-1}{2}}g$ is  smoothed %AG
			by $T$ from a $g\in \H$), then $a=\frac{2}{2+3}=0.4$ and the convergence rate is
			$\left[\log(N)/N\right]^{0.4}$; in other words the rate is approximately $1/N^{0.4}$. If $c$ takes its minimal value ($c=1$; $f_{\H}$ is less 
			smooth), then $a=\frac{1}{1+3}=\frac{1}{4}$ results in an approximate rate of $1/N^{0.25}$. 
			Alternatively, if we keep the total number of samples processed $t=lN=N^{a+1}$ fixed, $r(t) \approx 1/N^a = 1/t^{a/(a+1)} = 1/t^{1 - 1/(a+1)}$, i.e.,
			the convergence rate becomes larger for smoother regression problems (increasing $c$).
		  \item
			Last row: At this extreme, two terms dominate: the complexity of $\H$ according to $\rho_X$, and a term from the bound on $S_1$. Under this condition,
			although one can solve the matching criterion for $\lambda$, and it is possible to drive the individual terms of $r$ to zero, $l$ cannot be chosen large enough
			(within the analysed $l=N^a$ ($a>0$) scheme) to satisfy the $l\ge \lambda^{-\frac{1}{b}-1}$ constraint; thus convergence fails.
	      \end{compactenum}
	      \begin{table*}
                %\small
		\begin{center}
		\caption{Convergence conditions, convergence rates.
		Rows from top: $1-2$, $1-3$, $1-4$, $2-3$, $2-4$, $3-4^{th}$ terms are matched in $r(l,N,\lambda)$, the upper bound on the excess risk; see Eq.~\eqref{eq:r}. First column: convergence condition.
		      Second column: conditions for the dominance of the matched terms \emph{while} they also converge to zero. Third column: convergence rate of the excess risk.}\label{tab:conv-rates}
		\begin{tabular}{@{}ccc@{}} \toprule
		    Convergence condition & Dominance + convergence condition & Convergence rate\\\midrule
		    $\max\left(\frac{h}{(c+3)\min(2,b)},\frac{h(b+1)}{(c+3)b}\right)\le a$  & $\max\left(\frac{h\left(\frac{1}{b}+c\right)}{c+3},\frac{h(b+1)}{(c+3)b}\right)\le a$ & $\left[ \frac{\log(N)}{N}\right]^{\frac{hc}{c+3}}$\\
		    $\max\left( \frac{h}{6},\frac{h}{2(b+1)},\frac{h(b+1)}{2(2b+1)}\right) \le a < \frac{h}{2}$ & $\max\left(\frac{h}{6},\frac{h(b+1)}{2(2b+1)}\right)\le a < \min\left( \frac{h}{2}-\frac{h}{c+3}, \frac{\frac{h}{2}\left(\frac{1}{b}-1\right)}{\frac{1}{b}-2} \right)$ & $\frac{1}{N^{3a-\frac{h}{2}} \log^{\frac{h}{2}}(N)}$\\
		    $ \max\left(\frac{hb}{7b-2},\frac{h}{3b},\frac{h(b+1)}{4b}\right)\le a<h$ &  $\max\left(\frac{h(b-1)}{4b-2},\frac{h}{3b},\frac{h(b+1)}{4b}\right) \le a < \frac{h(bc+1)}{3b+bc}$ & $\frac{1}{N^{a+\frac{a-h}{3b-1}}\log^{\frac{h}{3b-1}}(N)}$\\
		    $a<\frac{h(c+1)}{6}$, $1>\frac{2(b+1)}{(c+1)b}$ &  never & never \vspace*{0.15cm}\\
		    $a<\frac{h(bc+1)}{3b}$, $1>\frac{b+1}{bc+1}$ & $a<\frac{h(bc+1)}{3b+bc}$, $1>\frac{b+1}{bc+1}$ & $\frac{1}{N^{\frac{abc}{bc+1}}}$\\
		     never &  never & never\\
		    \bottomrule
		\end{tabular}
		\end{center}
		\vspace*{-0.55cm}
	      \end{table*}

{\bf Proof of main theorem:} We present the main steps of the proof of our theorem; detailed derivations can be found in the supplementary material. 
Let us define $\b{x} = \{x_i\}_{i=1}^l$ and $\hat{\b{x}} = \{\{x_{i,n}\}_{n=1}^{N}\}_{i=1}^l$ as the `x-part' of $\b{z}$ and $\hat{\b{z}}$.
One can  express $f_{\b{z}}^{\lambda}$ \cite{caponnetto07optimal}, and similarly $f_{\hat{\b{z}}}^{\lambda}$ as
\begin{align}
    f_{\b{z}}^{\lambda} &= (T_{\b{x}}+\lambda)^{-1}g_{\b{z}}, \hspace*{0.65cm} T_{\b{x}} = \frac{1}{l}\sum_{i=1}^lT_{\mu_{x_i}},  \label{eq:f_zlambda}\\
    f_{\hat{\b{z}}}^{\lambda} &= (T_{\hat{\b{x}}}+\lambda)^{-1}g_{\hat{\b{z}}}, \hspace*{0.65cm} T_{\hat{\b{x}}} = \frac{1}{l}\sum_{i=1}^lT_{\mu_{\hat{x}_i}}, \label{eq:f_hatz^lambda}\\
    g_{\b{z}} &= \frac{1}{l}\sum_{i=1}^l K(\cdot,{\mu_{x_i}}) y_i, \hspace*{0.2cm} g_{\hat{\b{z}}} = \frac{1}{l}\sum_{i=1}^l K(\cdot,\mu_{\hat{x}_i}) y_i. \label{eq:gz,gzhat}
\end{align}
In Eqs.~\eqref{eq:f_zlambda}, \eqref{eq:f_hatz^lambda}, \eqref{eq:gz,gzhat}, $T_{\b{x}}, T_{\hat{\b{x}}}:\H\rightarrow \H$, $g_{\b{z}}, g_{\hat{\b{z}}}\in\H$. 

$\bullet$\hspace{1mm}\tb{Decomposition of the excess risk:}
	    We derive the upper bound for the excess risk
	    \begin{align}
		\hspace*{-0.2cm}\Eo\left[f^{\lambda}_{\hat{\b{z}}}\right] - \Eo\left[f_{\H}\right] \hspace*{-0.05cm}\le \hspace*{-0.05cm} 5 \left[S_{-1} + S_0 + \A(\lambda) + S_1 + S_2 \right]. \label{eq:E:5terms}
	    \end{align}
$\bullet$\hspace{1mm}\tb{It is sufficient to upper bound $S_{-1}$ and $S_0$:} \cite{caponnetto07optimal} has shown that $\forall \eta>0$ if
	      $l\ge \frac{2C_{\eta}B_K\N(\lambda)}{\lambda}$ and $\lambda \le \left\|T\right\|_{\L(\H)}$, then $\Pr(\bm{\Theta}(\lambda,\b{z})\le 1/2)\ge 1-\eta/3$, where
	      \begin{align*}
		  \bm{\Theta}(\lambda,\b{z}) = \|(T - T_{\b{x}}) (T+\lambda)^{-1}\|_{\L(\H)}
	      \end{align*}
	      and one can obtain upper bounds on $S_1$ and $S_2$ which hold with probability $1-\eta$.
      	      For $\A(\lambda)$ no probabilistic argument was needed.\vspace{1mm}\\
$\bullet$\hspace*{1mm}\tb{Probabilistic bounds on $\| g_{\hat{\b{z}}} - g_{\b{z}}\|_{\H}^2$, $\| T_{\b{x}}-T_{\hat{\b{x}}}\|_{\L(\H)}^2$, $\|\sqrt{T} (T_{\hat{\b{x}}}+\lambda)^{-1}\|_{\L(\H)}^2$, $\|f_{\b{z}}^{\lambda}\|_{\H}^2$:}
	By using the $\|Mu\|_{\H}\le \|M\|_{\L(\H)} \|u\|_{\H}$ $(M\in \L(\H), u\in \H)$ inequality, we bound $S_{-1}$ and $S_0$ as
	      \begin{align*}
		    S_{-1} &\le \| \sqrt{T} (T_{\hat{\b{x}}}+\lambda)^{-1}\|_{\L(\H)}^2 \|g_{\hat{\b{z}}} - g_{\b{z}} \|_{\H}^2,\\
		    S_0 &\le  \| \sqrt{T} (T_{\hat{\b{x}}}+\lambda)^{-1}\|_{\L(\H)}^2 \| T_{\b{x}}-T_{\hat{\b{x}}}\|_{\L(\H)}^2 \|f^{\lambda}_{\b{z}} \|^2_{\H}.
	      \end{align*}
	      For the terms on the r.h.s., we can derive the upper bounds [for $\alpha$ see Eq.~\eqref{eq:emp-mean-emb-conv-rate}]:
	      \begin{align*}
		  \left\|  g_{\hat{\b{z}}} - g_{\b{z}}  \right\|_{\H}^2  &\le  L^2 C^2  \frac{\left(1+\sqrt{\alpha}\right)^{2h} (2B_k)^h}{N^h},\\ 
		  \| \sqrt{T} (T_{\hat{\b{x}}}+\lambda)^{-1}\|_{\L(\H)}  &\le \frac{2}{\sqrt{\lambda}},\\
		  \left\|T_{\b{x}} - T_{\hat{\b{x}}}\right\|_{\L(\H)}^2 &\le  \frac{\left(1+\sqrt{\alpha}\right)^{2h} 2^{h+2}(B_k)^{h}B_K L^2}{N^h},\\
		  \left\|f_{\b{z}}^{\lambda}\right\|_{\H}^2 &\le \frac{C^2 B_K}{\lambda^2}.
	      \end{align*}
	      The bounds hold under the following conditions:
	      \begin{compactenum}
		  \item $\| g_{\hat{\b{z}}} - g_{\b{z}}\|_{\H}^2$: if the empirical mean embeddings are close to their population counterparts, i.e.,
			    \begin{align}
				\left\|\mu_{x_i} - \mu_{\hat{x}_i}\right\|_H &\le \frac{(1+\sqrt{\alpha})\sqrt{2B_k}}{\sqrt{N}} \label{eq:emp-mean-emb-conv-rate}
			  \end{align}
			  for $\forall i=1,\ldots,l$. This event has probability $1-le^{-\alpha}$ over all $l$ samples by a union bound.
		  \item $\| T_{\b{x}}-T_{\hat{\b{x}}}\|_{\L(\H)}^2$: \eqref{eq:emp-mean-emb-conv-rate} is assumed.
		  \item\label{eq:bounds} $\|\sqrt{T} (T_{\hat{\b{x}}}+\lambda)^{-1}\|_{\L(\H)}^2$: $\frac{\left(1+\sqrt{\alpha}\right)^{2} 2^{\frac{h+6}{h}}B_k (B_K)^{\frac{1}{h}} L^{\frac{2}{h}}}{(\lambda)^{\frac{2}{h}}}  \le N$,
			    \eqref{eq:emp-mean-emb-conv-rate},  and $\bm{\Theta}(\lambda,\b{z})\le \frac{1}{2}$.
		  \item $\|f_{\b{z}}^{\lambda}\|_{\H}^2$:  This upper bound always holds  (under the model assumptions).
	      \end{compactenum}
$\bullet$\hspace*{1mm}\tb{Union bound:} By applying an $\alpha=\log(l)+\delta$ reparameterization, and combining the received upper bounds with \cite{caponnetto07optimal}'s results for $S_1$ and $S_2$, the theorem follows with a union bound.

Finally, we note that 
\begin{compactitem}
  \item existing results were used at two points to simplify our analysis: bounding
	$S_1$, $S_2$, $\Theta(\lambda,\b{z})$ \cite{caponnetto07optimal} and $\left\|\mu_{x_i}-\mu_{\hat{x}_i}\right\|_H$ \cite{altun06unifying}.
 \item although the primary focus of our paper is clearly theoretical, we have provided some illustrative experiments in the supplementary material. These include 
	\begin{compactenum}
	  \item a comparison with the only alternative, theoretically justified distribution regression method \cite{poczos13distribution}\textsuperscript{\ref{footnote:RNDkitchensink}} on supervised entropy learning, where our approach gives better performance,
	  \item an experiment on aerosol prediction based on satellite images, where we perform as well as recent domain-specific, engineered methods \cite{wang12mixture} (which themselves beat state-of-the-art multiple instance learning alternatives).
	\end{compactenum}
\end{compactitem}

\vspace*{0.75cm}
\section{CONCLUSION}\label{sec:conclusions}
\vspace*{-0.1cm}
In this paper we established the learning theory of distribution regression under mild conditions, for
probability measures on separable, topological domains endowed with kernels.
We analysed an algorithmically  simple and parallelizable\footnote{Recently, \cite{zhang14divide} has constructed theoretically sound parallelization algorithms for kernel ridge regression.}
ridge regression scheme defined on the embeddings of the input distributions to a RKHS.
As a special case of our analysis, we proved the consistency of regression for set kernels \cite{haussler99convolution,gartner02multi} in the distribution-to-real regression 
setting (which was a 15-year-old open problem), and  for a recent kernel family \cite{christmann10universal}, which we have expanded upon (Table \ref{tab:K examples}).
To keep the presentation simple we focused on the quadratic loss ($\Eo$), bounded kernels ($k$, $K$), real-valued labels ($Y$), and mean embedding ($\mu$) based distribution 
regression with i.i.d.\ samples ($\{x_{i,n}\}_{n=1}^N$). In future work, we will relax these assumptions, and also consider deriving bounds with 
approximation error (capturing the richness of class $\H$ in the bounds).\footnote{The extension to separable Hilbert output spaces and the misspecicified case with approximation error are already available \cite{szabo14learning}.} Another exciting open question is whether (i) lower bounds on convergence can be proved, (ii)
optimal convergence rates can be derived,   (iii) one can obtain error bounds for non-point estimates.

\vspace*{-0.1cm}
\subsubsection*{Acknowledgements} This work was supported by the Gatsby Charitable Foundation, and by NSF grants IIS1247658 and IIS1250350. The work was carried out while Bharath K. Sriperumbudur was a research fellow in the Statistical Laboratory, Department of Pure Mathematics and Mathematical Statistics at the University of Cambridge, UK.
\bibliography{./BIB/MERR_short}
\bibliographystyle{unsrt}

%***********************************
%%supplementary material:
\newpage \onecolumn \appendix
\setdefaultleftmargin{2.5em}{2.2em}{1.87em}{1.7em}{1em}{1em}%paralist doc: default
\section{SUPPLEMENTARY MATERIAL}
This supplementary material contains (i) detailed proofs of the consistency of MERR (Section~\ref{appendix:proofs}), (ii) numerical illustrations (Section~\ref{sec: numerical experiments}).

\subsection{Proofs}\label{appendix:proofs}
\subsubsection{Proof of $k$: continuous, bounded $\Rightarrow$ $\mu$: $H$-measurable;  $\mu$: $H$-measurable, $X=\mu(\M^+_1(\X))\in \Bo(H)$ $\Rightarrow$  $\mu$: $X$-measurable $\Rightarrow \exists \rho$}\label{proof:rho:well-defined} 
Below we give sufficient conditions for the existence of probability measure $\rho$. We divide the proof into $3$ steps:

\begin{compactitem}
  \item \tb{$k$: continuous, bounded $\Rightarrow$ $\mu$: $H$-measurable}: The mapping 
	  $\mu: (\M^+_1(\X),\Bo(\tau_w))\rightarrow (H,\Bo(H))$ is measurable, iff the $L_h: (\M^+_1(\X),\Bo(\tau_w))\rightarrow (\R,\Bo(\R))$ map defined as $L_h(x) =\left<h,\mu_x\right>_H \left(=\int_{\X}h(u)\d x(u)\right)$ is measurable for $\forall h\in H$ \cite[Theorem IV.~22, page~116]{reed80methods}.
	  If $k$ is assumed to be continuous and bounded, these properties also hold for $\forall h\in H$ \cite[Lemma~4.23, page~124; Lemma~4.28, page~128]{steinwart08support}, i.e. $H=H(k)\subseteq C_b(\X)$. By the definition of 
	  the weak topology the $L_h$ functions are continuous (for $\forall h\in H$), which implies the required Borel measurability \cite[page~480]{steinwart08support} of $L_h$-s (for $\forall h\in H$).
  \item \tb{$\mu:$ $H$-measurable, $X=\mu(\M^+_1(\X)\in \Bo(H)$ $\Rightarrow$ $\mu:$ $X$-measurable}:
	Let $\tau$ denote the open sets on $H=H(k)$. Let $\left.\tau\right|_X=\{A\cap X: A \in \tau\}$ be the subspace topology on $X$, and let 
	$\left.\Bo(H)\right|_X=\{A\cap X: A\in \Bo(H)\}$ be the subspace $\sigma$-algebra on $X$. Since $\Bo\left(\left.\tau\right|_X\right) = \left.\Bo(H)\right|_X \subseteq \Bo(H)$ 
	(the containing relation follows from the $X\in\Bo(H)$ condition), and $\left.\Bo(H)\right|_X=\{A\in  \Bo(H): A\subseteq X\}$, the 
	measurability of $\mu: (\M^+_1(\X),\Bo(\tau_w))\rightarrow (H,\Bo(H))$ implies the measurability of $\mu: (\M^+_1(\X),\Bo(\tau_w))\rightarrow (X,\left.\Bo(H)\right|_X)$.
  \item \tb{$\mu:$ $X$-measurable  $\Rightarrow$ $\exists \rho$}: Let us consider the
	\begin{align}
	    g:(\M^+_1(\X)\times \R, \Bo(\tau_w) \otimes \Bo(\R)) \rightarrow (X\times \R,\left.\Bo(H)\right|_X \otimes \Bo(\R)) \label{eq:g}
	\end{align}
	$g(x,y) = [g_1(x,y);g_2(x,y)] = [\mu_x;y]$ mapping. If $g$ is a \emph{measurable} function, then it defines $\rho$, a probability measure on 
	$(X\times  \R,\left.\Bo(H)\right|_X \otimes \Bo(\R))$ 
	by looking at $g$ as a random variable taking values in $X\times \R$:
	\begin{align}
	  \rho(C) &:= \M\left(g^{-1}(C)\right), \quad (C\in \left.\Bo(H)\right|_X\otimes \Bo(\R)). \label{eq:rho-def}
	\end{align}
	Function $g$ in Eq.~\eqref{eq:g} is measurable iff its coordinate functions, $g_1$ and $g_2$ are both measurable functions \cite[Proposition~3.2, page 201]{nagyXreal}.
	  Thus, we need for $\forall A\in \left.\Bo(H)\right|_X$, $\forall B\in \Bo(\R)$
	  \begin{align}
	  \Bo(\tau_w)\otimes \Bo(\R)&\ni g_1^{-1}(A) = \{(x,y): g_1(x,y)=\mu_x\in A\}=\mu^{-1}(A)\times \R,\label{eq:g1}\\
	  \Bo(\tau_w)\otimes \Bo(\R)&\ni g_2^{-1}(B) = \{(x,y): g_2(x,y)=y\in B\} = \M^+_1(\X)\times B.\label{eq:g2}
	  \end{align}
	  According to Eqs.~\eqref{eq:g1}-\eqref{eq:g2}, the measurability of $g$ follows from the $X$-measurability of 
	  $\mu: (\M^+_1(\X),\Bo(\tau_w))\rightarrow (X,\left.\Bo(H)\right|_X)$, which is guaranteed by our conditions. 
\end{compactitem}

\subsubsection{Proof of $\Psi$: H{\"o}lder continuous $\Rightarrow$ $K$: measurable}
  \cite{caponnetto07optimal}'s original assumption that $(\mu_a,\mu_b)\in X \times X \mapsto K(\mu_a,\mu_b)\in\R$ is measureable follows from the required H{\"o}lder continuity [see Eq.~\eqref{eq:K:Lip}]
    since (i) the continuity of $\Psi$ is equivalent to that of $K$, (ii) a continuous map between topological spaces is Borel measurable \cite[Lemma~4.29 on page 128; page 480]{steinwart08support}.

\subsubsection{Proof of $K$: linear $\Rightarrow$ $\Psi$: H{\"o}lder continuous with $L=1$, $h=1$}

	In case of a linear $K$ kernel
	      $K(\mu_a,\mu_b)=\left<\mu_a,\mu_b\right>_H \quad (\mu_a,\mu_b\in X)$,
	      by the bilinearity of $\left<\cdot,\cdot\right>_H$ and $\left\|\left<\cdot,a\right>_H\right\|_{\H}^2 = \left\|a\right\|_{H}^2$,
	      we get that
              $\left\|K(\cdot,\mu_a) - K(\cdot,\mu_b)\right\|_{\H} =\left\|\left<\cdot,\mu_a\right>_H - \left<\cdot,\mu_b\right>_H\right\|_{\H}
		  =\left\|\left<\cdot,\mu_a-\mu_b\right>_H\right\|_{\H} = \left\|\mu_a-\mu_b\right\|_H$.
	      In other words, H{\"o}lder continuity holds with $L=1$, $h=1$; $K$ is Lipschitz continuous ($h=1$).

\subsubsection{Proof of $\X$: compact metric, $\mu$: continuous $\Rightarrow$ $X=\mu\left(\M^+_1(\X)\right)$: compact metric}\label{proof:X compact}
    Let us suppose that $\X=(\X,d)$ is a compact metric space. This implies that $\M^+_1(\X)$ is also a compact metric space by Theorem~6.4 in \cite{parthasarathy67probability} (page~55). 
     The continuous ($\mu$) image of a compact set is compact 
    (see page~478 in \cite{steinwart08support}), thus $X=\mu\left(\M^+_1(\X)\right)\subseteq H$ is compact metric.

\subsubsection{Proof of the Kernel Examples on $X = \mu(\M^+_1(\X))$}\label{proof:Kexamples}
Below we prove for the $K:X\times X\rightarrow \R$ functions in Table~\ref{tab:K examples} that they are kernels on mean embedded distributions. 

We need some definitions and lemmas. $\Z,\Z^+,\R^+,\R^{\ge 0}$ denotes the set of integers, positive integers, positive real numbers and non-negative real numbers, respectively.
	      \begin{definition}
		    Let $X$ be a non-empty set. A $K: X \times X \rightarrow \R$ function is called 
		    \begin{compactitem}
		      \item \emph{positive definite} (pd; also referred to as kernel) on $X$, if it is 
			  \begin{compactenum}
			      \item symmetric [$K(a,b)=K(b,a)$, $\forall a,b\in X$], and 
			      \item $\sum_{i,j=1}^nc_ic_j K(a_i,a_j)\ge 0$ for all $n\in\Z^+$, $\{a_1,\ldots,a_n\}\subseteq X^n$, $\b{c}=[c_1;\ldots; c_n]\in\R^n$.
			  \end{compactenum}
		      \item \emph{negative definite} (nd; sometimes $-K$ is called conditionally positive definite) on $X$, if it is 
			  \begin{compactenum}
			      \item symmetric, and
			      \item $\sum_{i,j=1}^nc_ic_j K(a_i,a_j)\le 0$ for all $n\in\Z^+$, $\{a_1,\ldots,a_n\}\subseteq X^n$, $\b{c}=[c_1;\ldots;c_n]\in\R^n$, where $\sum_{j=1}^n c_j=0$.
			  \end{compactenum}
		    \end{compactitem}
	      \end{definition}

	      We will use the following properties of positive/negative definite functions:
	      \begin{compactenum}
		\item\label{prop:exp} $K$ is nd $\Leftrightarrow$ $e^{-tK}$ is pd for all $t>0$; see Chapter~3 in \cite{berg84harmonic}.
		\item\label{prop:t+inverse} $K:X\times X \rightarrow \R^{\ge 0}$ is nd $\Leftrightarrow$ $\frac{1}{t+K}$ is pd for all $t>0$; see Chapter~3 in \cite{berg84harmonic}.
		\item\label{prop:exponent} If $K$ is nd and non-negative on the diagonal ($K(x,x)\ge 0$, $\forall x\in X$), then $K^{\alpha}$ is nd for all $\alpha\in [0,1]$; see Chapter~3 in \cite{berg84harmonic}.
		\item\label{prop:ip:pd} $K(x,y)=\left<x,y\right>_X$ is pd, where $X$ is a Hilbert space (since the pd property is equivalent to being a kernel).
		\item\label{prop:norm:nd} $K(x,y)=\left\|x-y\right\|_X^2$ is nd, where $X$ is a Hilbert space; see Chapter~3 in \cite{berg84harmonic}.
		\item\label{prop:ip+const} If $K$ is nd, $K+d$ ($d\in \R$) is also nd. Proof: (i) $K(x,y)+d=K(y,x)+d$ holds by the symmetry of $K$, (ii) $\sum_{i,j=1}^nc_ic_j \left[K(a_i,a_j)+d\right] = \sum_{i,j=1}^nc_ic_j K(a_i,a_j) + \sum_{i=1}^nc_i\sum_{j=1}^nc_j d = \sum_{i,j=1}^nc_ic_j K(a_i,a_j) + \sum_{i=1}^nc_id\sum_{j=1}^nc_j = \sum_{i,j=1}^nc_ic_j K(a_i,a_j)+0\le 0$, where we used that $\sum_{j=1}^nc_j=0$ and $K$ is nd.
		\item\label{prop:restriction} If $K$ is pd (nd) on $X$, then it is pd (nd) on $X'\subseteq X$ as well. Proof: less constraints have to be satisfied for $X'\subseteq X$.
		\item\label{prop:posmult} If $K$ is pd (nd) on $X$, then $sK$ ($s\in\R^+$) is also pd (nd). Proof: multiplication by a positive constant does not affect the sign of $\sum_{i,j=1}^nc_ic_j K(a_i,a_j)$. 
		\item\label{prop:inverse} If $K$ is nd on $X$ and $K(x,y)>0$ $\forall x,y\in X$, then $\frac{1}{K}$ is pd; see Chapter~3 in \cite{berg84harmonic}.
		\item\label{prop:hol-composition} If $K$ is pd on $X$, and $h(u)=\sum_{n=0}^{\infty}a_nu^n$ with $a_n\ge 0$, then $h\circ K$ is pd; see Chapter~3 in \cite{berg84harmonic}.
	      \end{compactenum}

	      Making use of these properties one can prove the kernel property of the $K$-s in Table~\ref{tab:K examples} (see also Table~\ref{tab:K examples:more generally}) as follows.
	      All the $K$-s are functions of $\left\|\mu_a-\mu_b\right\|_H$, $\left\|\mu_a-\mu_b\right\|_H = \left\|\mu_b-\mu_a\right\|_H$, hence $K$-s are symmetric.

	      $K(x,y)=\left\|x-y\right\|_H^2$ is nd on $H=H(k)$ (Prop.~\ref{prop:norm:nd}), thus $K(x,y)=\left\|x-y\right\|_H^2$ is nd on $X=\mu\left(\M^+_1(\X)\right)\subseteq H(k)$ (Prop.~\ref{prop:restriction}). 
	      Consequently, $K(x,y)=\left\|x-y\right\|_H^d$ is nd on $X$, where $d\in[0,2]$ ($K(x,x)=0\ge 0$, Prop.~\ref{prop:exponent}). 
	      \begin{compactitem}
		\item Hence, $K(x,y)=e^{-t\left\|x-y\right\|_H^d}$ is pd, where $t>0$, $d\in[0,2]$ (Prop.~\ref{prop:exp}). By the $(t,d)=\left(\frac{1}{2\theta^2},2\right)$ and $(t,d)=\left(\frac{1}{2\theta^2},1\right)$ choices, we get that $K_G$ and $K_e$ are kernels.
		\item Using Prop.~\ref{prop:t+inverse} ($\left\|x-y\right\|_H^d\ge 0$), one obtains that $K(x,y)=\frac{1}{t+\left\|x-y\right\|_H^d}$ is pd on $X$, where $t>0$, $d\in[0,2]$. By the $(t,d)=(1,\le 2)$ choice the kernel property of $K_t$ follows.
		\item Thus, $K(x,y)=s \left\|x-y\right\|_H^d$ is nd on $X$, where $s>0$, $d\in [0,2]$ (Prop.~\ref{prop:posmult}). Consequently, $K(x,y)=\frac{1}{t+s \left\|x-y\right\|_H^d}$ is pd on $X$, where $s>0$, $d\in [0,2]$, $t>0$ (Prop.~\ref{prop:t+inverse}). 
		      By the $(d,t,s)=(2,1,\frac{1}{\theta^2})$, we have that $K_C$ is kernel.
		\item Hence,  $K(x,y)=\left\|x-y\right\|_H^d+e$ is nd on $X$, where $d\in[0,2]$, $e\in\R^+$ (Prop.~\ref{prop:ip+const}). Thus, $K(x,y)=\left(\left\|x-y\right\|_H^d+e\right)^f$ is nd on $X$, where $d\in[0,2]$, $e\in\R^+$, $f\in (0,1]$ ($\left\|x-y\right\|_H^d+e\ge 0$, Prop.~\ref{prop:exponent}). 
		      Consequently, $K(x,y)=\frac{1}{\left(\left\|x-y\right\|_H^d+e\right)^f}$ is pd on $X$, where $d\in[0,2]$, $e\in\R^+$, $f\in (0,1]$ ($\left(\left\|x-y\right\|_H^d+e\right)^f>0$; Prop.~\ref{prop:inverse}); with the $(d,e,f)=\left(2,\theta^2,\frac{1}{2}\right)$ choice, one obtains that $K_i$ is a kernel.
	      \end{compactitem}
	    \begin{table}
	      \begin{center}
	      \caption{Nonlinear kernels on mean embedded distributions.}
	      \label{tab:K examples:more generally}
	      \begin{tabular}{@{}r@{\hspace{0.1cm}}c@{\hspace{0.1cm}}ll@{}}
		\toprule
		\multicolumn{3}{c}{Kernel ($K$)} &  Parameter(s)\\ \midrule
		$K(\mu_a,\mu_b)$ & $=$ & $e^{-t\left\|\mu_a-\mu_b\right\|_H^d}$ & $t>0$, $d\in[0,2]$ \\
		$K(\mu_a,\mu_b)$ & $=$ & $\frac{1}{t+\left\|\mu_a-\mu_b\right\|_H^d}$ & $t>0$, $d\in[0,2]$ \\
		$K(\mu_a,\mu_b)$ & $=$ & $\frac{1}{t+s \left\|\mu_a-\mu_b\right\|_H^d}$ & $s>0$, $d\in [0,2]$, $t>0$\\ 
		$K(\mu_a,\mu_b)$ & $=$ & $\frac{1}{\left(\left\|\mu_a-\mu_b\right\|_H^d+e\right)^f}$ &  $d\in[0,2]$, $e\in\R^+$\\\bottomrule
	      \end{tabular}
	      \end{center}
	    \end{table}

\subsubsection{Proof of ``Conditions of Proof~\ref{proof:X compact} and Proof~\ref{proof:Kexamples}'' $\Rightarrow$ $\Psi$-s of $K$-s in Proof~\ref{proof:Kexamples}: H{\"o}lder continuous} 
We tackle the problem more generally: 
      \begin{compactenum}
	  \item 
	      we give sufficient conditions for $K$ kernels of the form
	    \begin{align}
		K(\mu_a,\mu_b) &= \bar{K}\left(\left\|\mu_a-\mu_b\right\|_H\right), \label{eq:radialK}
	    \end{align}
	      i.e., for radial kernels to  have H{\"o}lder continuous canonical feature map ($\Psi(\mu_c)=K(\cdot,\mu_c)$):
	      $\exists L>0$, $h\in(0,1]$ such that $\left\|K(\cdot,\mu_a) - K(\cdot,\mu_b)\right\|_{\H} \le L \left\|\mu_a - \mu_b\right\|_H^h$.
	  \item Then we show that these sufficient conditions are satisfied for the $K$ kernels listed in Table~\ref{tab:K examples}. 
      \end{compactenum}
      Let us first note that  $K$ is bounded. Indeed, since $\Psi$ is H{\"o}lder continuous, specially it is continuous. Hence using Lemma 4.29 in \cite{steinwart08support} (page 128), the  
     \begin{align*}
	  K_0:\mu_a\in X \rightarrow K(\mu_a,\mu_a)\in \R
     \end{align*}
     mapping is continuous. As we have already seen (Section~\ref{proof:X compact}) $X$ is compact. The continuous ($K_0$) image of a compact set ($X$), i.e., the $\{K(\mu_a,\mu_a): \mu_a\in X\}\subseteq \R$ set is compact, specially it is bounded above.

      \begin{compactenum}
	  \item Sufficient conditions:       
		  Now, we present sufficient conditions for the assumed H{\"o}lder continuity
		    \begin{align}
			  \left\|K(\cdot,\mu_a) - K(\cdot,\mu_b)\right\|_{\H} &\le L \left\|\mu_a - \mu_b\right\|_H^h. \label{eq:K:Holder}
		    \end{align}
		  Using $\left\|u\right\|_{\H}^2=\left<u,u\right>_{\H}$, the bilinearity of $\left<\cdot,\cdot\right>_{\H}$, 
		  the reproducing property of $K$ and Eq.~\eqref{eq:radialK}, we get
		  \begin{align*}
			\left\|K(\cdot,\mu_a) - K(\cdot,\mu_b)\right\|_{\H}^2 &= \left<K(\cdot,\mu_a) - K(\cdot,\mu_b),K(\cdot,\mu_a) - K(\cdot,\mu_b)\right>_{\H}\\
										&= K(\mu_a,\mu_a) + K(\mu_b,\mu_b) -2 K(\mu_a,\mu_b)
										= 2\bar{K}(0) - 2 \bar{K}(\left\|\mu_a-\mu_b\right\|_H)\\
										& = 2 \left[\bar{K}(0) -\bar{K}\left(\left\|\mu_a-\mu_b\right\|_H\right)\right].
		  \end{align*}
		  Hence, the H{\"o}lder continuity of $K$ is equivalent to the existence of an $L'\left(=\frac{L^2}{2}\right)>0$ such that
		  \begin{align*}
		      \bar{K}(0) -\bar{K}\left(\left\|\mu_a-\mu_b\right\|_H\right) & \le L' \left\|\mu_a-\mu_b\right\|_H^{2h}.
		  \end{align*}
		  Since for $\mu_a=\mu_b$ both sides are equal to $0$, this requirement is equivalent to
		  \begin{align*}
		      u(\mu_a,\mu_b): = \frac{\bar{K}(0) -\bar{K}\left(\left\|\mu_a-\mu_b\right\|_H\right)}{\left\|\mu_a-\mu_b\right\|_H^{2h}} \le L', \quad(\mu_a\ne \mu_b)
		  \end{align*}
		    i.e., that the $u:X\times X\rightarrow \R$ function is bounded above. Function $u$ is the composition ($u=u_2\circ u_1$) of the mappings:
		  \begin{align}
		      u_1&: X\times X \rightarrow \R^{\ge 0}, \quad u_1(\mu_a,\mu_b) = \left\|\mu_a-\mu_b\right\|_H,\nonumber\\
		      u_2&: \R^{\ge 0} \rightarrow \R, \quad u_2(v) = \frac{\bar{K}(0)-\bar{K}(v)}{v^{2h}}. \label{eq:u_2}
		  \end{align}
		  Here, $u_1$ is continuous. Let us \emph{suppose} for $u_2$ that 
   	         \begin{compactenum}
				\item (i) $\exists h\in (0,1]$ such that $\lim_{v\rightarrow 0+}u_2(v)$ exists, and		 
				\item $u_2$ is continuous.
		  \end{compactenum}
		  In this case, since the composition of continuous functions is continuous (see page 85 in \cite{kelley75general}), $u$ is continuous. As we have seen (Section~\ref{proof:X compact}), $X$ is compact. The product of compact sets ($X\times X$) is compact by the Tychonoff theorem (see page 143 in  \cite{kelley75general}).
		  Finally, since the continuous ($u$) image of a compact set ($X\times X$), i.e. $\{u(\mu_a,\mu_b): (\mu_a,\mu_b)\in X\times X\}\subseteq\R$  is compact (Theorem~8 in \cite{kelley75general}, page 141), we get that $u$ is bounded, specially bounded above.

		  To sum up, we have proved that if
		  \begin{compactenum}
		      \item $K$ is radial [see Eq.~\eqref{eq:radialK}],
		      \item $u_2$ [Eq.~\eqref{eq:u_2}] is (i) continuous and (ii) $\exists h\in (0,1]$ such that $\lim_{v\rightarrow 0+}u_2(v)$ exists,
		  \end{compactenum}
		   then the H{\"o}lder property [Eq.~\eqref{eq:K:Holder}] holds for $K$ with exponent $h$. In other words, the H{\"o}lder property of a kernel $K$ on mean embedded distributions can be simply guaranteed by the appropriate behavior of $\bar{K}$ at zero.
	  \item Verification of the sufficient conditions: In the sequel we show that these conditions hold for the $u_2$ functions of the $K$ kernels in Table~\ref{tab:K examples}. In the examples
	      \begin{align*}
		    \bar{K}_G(v) &= e^{-\frac{v^2}{2\theta^2}}, &  \bar{K}_e(v) &= e^{-\frac{v}{2\theta^2}}, &  \bar{K}_C(v) &= \frac{1}{1+\frac{v^2}{\theta^2}}, &
		    \bar{K}_t(v) &= \frac{1}{1+v^{\theta}}, & \bar{K}_i(v) &= \frac{1}{\sqrt{v^2+\theta^2}}.
	      \end{align*}
	      The corresponding $u_2$ functions are
	      \begin{align*}
		    u_{2G}(v) &= \frac{1-e^{-\frac{v^2}{2\theta^2}}}{v^{2h}}, & 
		    u_{2e}(v) &= \frac{1-e^{-\frac{v}{2\theta^2}}}{v^{2h}}, &
		    u_{2C}(v) &= \frac{1-\frac{1}{1+\frac{v^2}{\theta^2}}}{v^{2h}}, &
		    u_{2t}(v) &= \frac{1-\frac{1}{1+v^{\theta}}}{v^{2h}}, &
		    u_{2i}(v) &= \frac{\frac{1}{\theta}-\frac{1}{\sqrt{v^2+\theta^2}}}{v^{2h}}.
	      \end{align*}
	      The limit requirements at zero complementing the continuity of $u_2$-s are satisfied:
	      \begin{compactitem}
		\item $u_{2G}$: In this case
		    \begin{align*}
			  \lim_{v\rightarrow 0+}u_{2G}(v) &= \lim_{v\rightarrow 0+} \frac{1-e^{-\frac{v^2}{2\theta^2}}}{v^{2}} =  \lim_{v\rightarrow 0+} \frac{1-e^{-\frac{v}{2\theta^2}}}{v} = \lim_{v\rightarrow 0+} \frac{\frac{1}{2\theta^2}e^{-\frac{v}{2\theta^2}}}{1} = \frac{1}{2\theta^2},	  
		    \end{align*}
		    where we applied a $v^2$ substitution and the L'Hopital rule; $h=1$.
		\item $u_{2e}$:
		    \begin{align*}
			    \lim_{v\rightarrow 0+}u_{2e}(v) &= \lim_{v\rightarrow 0+}\frac{1-e^{-\frac{v}{2\theta^2}}}{v^{2h}} = \lim_{v\rightarrow 0+}\frac{\frac{1}{2\theta^2}e^{-\frac{v}{2\theta^2}}}{2hv^{2h-1}}=\frac{1}{2\theta^2},
		    \end{align*}
		    where we applied the L'Hopital rule and chose $h=\frac{1}{2}$, the largest $h$ from the $2h-1\le 0$ convergence domain.
		\item $u_{2C}$:
		    \begin{align*}
			  u_{2C}(v) &= \frac{1-\frac{1}{1+\frac{v^2}{\theta^2}}}{v^{2h}} = \frac{1-\frac{\theta^2}{\theta^2+v^2}}{v^{2h}} = \frac{\frac{v^2}{\theta^2+v^2}}{v^{2h}} = \frac{v^{2-2h}}{\theta^2+v^2\xrightarrow{v\rightarrow 0+}\theta^2},
		    \end{align*}
		    we chose $h=1$, the largest value from the convergence domain ($2-2h\ge 0 \Rightarrow 1 \ge h$).
		\item $u_{2t}$:
		    \begin{align*}
			  u_{2t}(v) &= \frac{1-\frac{1}{1+v^{\theta}}}{v^{2h}} = \frac{v^{\theta-2h}}{1+v^{\theta}\xrightarrow{v\rightarrow 0} 1},  
		    \end{align*}
		    thus we can have $h=\frac{\theta}{2}$, the largest element of the convergence domain 
		    ($\theta-2h\ge 0 \Leftrightarrow \frac{\theta}{2}\ge h$). Here we require $\theta\le 2$ in order to guarantee that $h=\frac{\theta}{2}\le 1$.
		\item $u_{2i}$: Let $g$ denote the nominator of $u_{2i}$
		    \begin{align*}
			g(v) &=\frac{1}{\theta}-\frac{1}{\sqrt{v^2+\theta^2}}=\frac{1}{\theta}-\left[g(0)+g'(0)v + \frac{g''(0)}{2}v^2+\ldots\right]\\
			     &= \frac{1}{\theta}-\left[\frac{1}{\theta} + \left(\left.-\frac{1}{2}\frac{1}{(v^2+\theta)^{\frac{3}{2}}}2v\right)\right|_{v=0} v +\frac{g''(0)}{2}v^2+\ldots\right]
			     = -v^2\left[\frac{g''(0)}{2!}+\frac{g^{(3)}(0)}{3!}v+\ldots\right].
		    \end{align*}
		    Hence,
		    \begin{align*}
			     \lim_{v\rightarrow 0+}u_{2i}(v) &= \lim_{v\rightarrow 0+}\frac{\frac{1}{\theta}-\frac{1}{\sqrt{v^2+\theta^2}}}{v^{2}} = \lim_{v\rightarrow 0+}\frac{-v^2\left[\frac{g''(0)}{2!}+\frac{g^{(3)}(0)}{3!}v+\ldots\right]}{v^2} = -\frac{g''(0)}{2},
		    \end{align*}
		    i.e., $h$ can be chosen to be $1$ ($h=1$).
	      \end{compactitem}
      \end{compactenum}

\subsubsection{Proof of $\left\|\sum_{i=1}^nf_i\right\|^2 \le n \sum_{i=1}^n \left\|f_i\right\|^2$} \label{lemma:normeq}
 In a normed space $(N,\left\|\cdot\right\|)$
  \begin{align}
      \left\|\sum_{i=1}^nf_i\right\|^2 &\le n \sum_{i=1}^n \left\|f_i\right\|^2, \label{eq:norm-eq}
  \end{align}
  where $f_i\in N$ ($i=1,\ldots,n$). 

Indeed the statement holds since $\left\|\sum_{i=1}^nf_i\right\|^2 \le \left(\sum_{i=1}^n\left\|f_i\right\|\right)^2 \le n \sum_{i=1}^n \left\|f_i\right\|^2$,
where we applied the triangle inequality, and a consequence that the arithmetic mean is smaller or equal than the squared mean (special case of the generalized mean inequality) with $a_i =\left\|f_i\right\|\ge 0$. Particularly,
$\frac{\sum_{i=1}^na_i}{n} \le \sqrt{\frac{\sum_{i=1}^n(a_i)^2}{n}} \Rightarrow \left(\sum_{i=1}^na_i\right)^2 \le n \sum_{i=1}^n(a_i)^2.$

\subsubsection{Proof of the Decomposition of the Excess Risk}
    It is known \cite{caponnetto07optimal} that $\Eo[f] - \Eo[f_{\H}] = \| \sqrt{T}(f-f_{\H})\|_{\H}^2\quad (\forall f\in \H)$.
    Applying this identity with $f=f_{\hat{\b{z}}}^{\lambda}\in\H$ and a telescopic trick, we get
    \begin{align}
	  \Eo\left[f^{\lambda}_{\hat{\b{z}}}\right] - \Eo\left[f_{\H}\right] = \left\|\sqrt{T}\left(f^{\lambda}_{\hat{\b{z}}} - f_{\H}\right)\right\|_{\H}^2
	  &= \left\|\sqrt{T}\left[\left(f^{\lambda}_{\hat{\b{z}}} - f^{\lambda}_{\b{z}}\right) + \left(f^{\lambda}_{\b{z}}- f^{\lambda}\right) + \left(f^{\lambda} - f_{\H}\right)\right]\right\|_{\H}^2\label{eq:target}. 
    \end{align}
    By Eqs.~\eqref{eq:f_zlambda}, \eqref{eq:f_hatz^lambda}, 
    and the operator identity $A^{-1}-B^{-1} = A^{-1}(B-A)B^{-1}$
    one obtains for the first term in Eq.~\eqref{eq:target}
    \begin{align*}
	  f^{\lambda}_{\hat{\b{z}}} - f^{\lambda}_{\b{z}}  &= (T_{\hat{\b{x}}}+\lambda)^{-1}g_{\hat{\b{z}}} - (T_{\b{x}}+\lambda)^{-1}g_{\b{z}}
	  = (T_{\hat{\b{x}}}+\lambda)^{-1}(g_{\hat{\b{z}}} - g_{\b{z}}) + (T_{\hat{\b{x}}}+\lambda)^{-1}g_{\b{z}} - (T_{\b{x}}+\lambda)^{-1}g_{\b{z}}\\
	&= (T_{\hat{\b{x}}}+\lambda)^{-1}(g_{\hat{\b{z}}} - g_{\b{z}}) + \left[(T_{\hat{\b{x}}}+\lambda)^{-1} - (T_{\b{x}}+\lambda)^{-1}\right]g_{\b{z}}\\
	&= (T_{\hat{\b{x}}}+\lambda)^{-1}(g_{\hat{\b{z}}} - g_{\b{z}}) +\left[(T_{\hat{\b{x}}}+\lambda)^{-1}(T_{\b{x}}-T_{\hat{\b{x}}}) (T_{\b{x}}+\lambda)^{-1}\right]g_{\b{z}}\\
	&= (T_{\hat{\b{x}}}+\lambda)^{-1}\left[(g_{\hat{\b{z}}} - g_{\b{z}}) + (T_{\b{x}}-T_{\hat{\b{x}}}) (T_{\b{x}}+\lambda)^{-1}g_{\b{z}} \right]
	= (T_{\hat{\b{x}}}+\lambda)^{-1}\left[(g_{\hat{\b{z}}} - g_{\b{z}}) + (T_{\b{x}}-T_{\hat{\b{x}}}) f_{\b{z}}^{\lambda}\right].
    \end{align*}
    Thus, we can rewrite the first term in \eqref{eq:target} as 
    \begin{align*}
    \sqrt{T}\left(f^{\lambda}_{\hat{\b{z}}} - f^{\lambda}_{\b{z}}\right) &=: f_{-1}+f_0,& f_{-1} &= \sqrt{T} (T_{\hat{\b{x}}}+\lambda)^{-1}(g_{\hat{\b{z}}} - g_{\b{z}}),&
	  f_{0}  &= \sqrt{T} (T_{\hat{\b{x}}}+\lambda)^{-1}(T_{\b{x}}-T_{\hat{\b{x}}}) f^{\lambda}_{\b{z}}.
    \end{align*}
    The second term in \eqref{eq:target} can be decomposed \cite{caponnetto07optimal}  as
      \begin{align*}
	    \sqrt{T}\left[\left(f^{\lambda}_{\b{z}}- f^{\lambda}\right) + \left(f^{\lambda} - f_{\H}\right)\right]
	    & = \sqrt{T}\left[(T_{\b{x}}+\lambda)^{-1}(g_{\b{z}}-T_{\b{x}}f_{\H})+
	    (T_{\b{x}}+\lambda)^{-1}(T-T_{\b{x}})(f^{\lambda}-f_{\H}) + (f^{\lambda}-f_{\H})\right]\\
	    &=:  f_1 +f_2 + f_3,
    \end{align*}
      where 
      \begin{align*}
	f_1 & = \sqrt{T}(T_{\b{x}}+\lambda)^{-1}(g_{\b{z}}-T_{\b{x}}f_{\H}), &
	f_2 & = \sqrt{T}(T_{\b{x}}+\lambda)^{-1}(T-T_{\b{x}})(f^{\lambda}-f_{\H}),& 
	f_3 &= \sqrt{T}(f^{\lambda}-f_{\H}).  
      \end{align*}
      Using these $f_i$ notations, \eqref{eq:target} can be upper bounded as
		  \begin{align}
			\Eo\left[f^{\lambda}_{\hat{\b{z}}}\right] - \Eo\left[f_{\H}\right] &= \left\|\sum_{i=-1}^3f_{i}\right\|_{\H}^2 \le 5 \sum_{i=-1}^3\left\|f_i\right\|_{\H}^2, \label{eq:telescopic}
		  \end{align}
		  exploiting  Section~\ref{lemma:normeq} ($\left\|\cdot\right\|^2=\left\|\cdot\right\|^2_{\H}$, $n=5$).
	  Consequently, introducing the 
    \begin{align*}
      S_{-1} &=S_{-1}(\lambda,\b{z},\hat{\b{z}}) = \|f_{-1}\|_{\H}^2, &
      S_0 &=S_0(\lambda,\b{z},\hat{\b{z}}) = \|f_{0}\|_{\H}^2,\\
      S_1 &=S_1(\lambda,\b{z}) = \|f_{1}\|_{\H}^2,&
      S_2&=S_2(\lambda,\b{z}) = \|f_{2}\|_{\H}^2,&
      \A(\lambda) &= \|f_{3}\|_{\H}^2,
    \end{align*}
    notations (for $\A(\lambda)$ see also Theorem~\ref{theo1}), 
    \eqref{eq:telescopic} can be rewritten as	
    \begin{align}
      \Eo\left[f^{\lambda}_{\hat{\b{z}}}\right] - \Eo\left[f_{\H}\right] &\le 5 \left[S_{-1} + S_0 + \A(\lambda) + S_1 + S_2 \right]. \label{eq:5-bound}
    \end{align}

\subsubsection{Proof of the Upper Bounding Terms of $S_{-1}$ and $S_0$}
	      Using the
	      \begin{align}
		  \|Mu\|_{\H}\le \|M\|_{\L(\H)} \|u\|_{\H}\quad(M\in \L(\H), u\in \H), \label{eq:opnorm-fit}
	      \end{align}
	      relation, we get
		  \begin{align*}
		      S_{-1} &\le \left\| \sqrt{T} (T_{\hat{\b{x}}}+\lambda)^{-1}\right\|_{\L(\H)}^2 \left\|g_{\hat{\b{z}}} - g_{\b{z}} \right\|_{\H}^2,\\
		      S_0 &\le \left\| \sqrt{T} (T_{\hat{\b{x}}}+\lambda)^{-1}\right\|_{\L(\H)}^2 \left\| (T_{\b{x}}-T_{\hat{\b{x}}}) f^{\lambda}_{\b{z}} \right\|^2_{\H}
			  \le  \left\| \sqrt{T} (T_{\hat{\b{x}}}+\lambda)^{-1}\right\|_{\L(\H)}^2 \left\| T_{\b{x}}-T_{\hat{\b{x}}}\right\|_{\L(\H)}^2 \left\|f^{\lambda}_{\b{z}} \right\|^2_{\H}.
		  \end{align*}

\subsubsection{Proof of the Convergence Rate of the Empirical Mean Embedding}  \label{lemma:ME:convergence}
The statement we prove is as follows.\cite{altun06unifying}\footnote{In the original result a factor of $2$ is missing from the denominator in the exponential function; we correct the proof here.} 
 
Let $\mu_x = \int_{\X}k(\cdot,u)\d x(u)$ denote the mean embedding of  distribution $x\in \M^+_1\left(\X\right)$ to the $H = H(k)$ RKHS determined by kernel $k$ ($\mu_x\in H$), which is assumed to be bounded $k(u,u)\le B_k$ ($\forall u\in \X$). 
Let us given $N$ i.i.d.\ samples from distribution $x$: $x_1$, ..., $x_N$. Let $\mu_{\hat{x}}=\frac{1}{N}\sum_{n=1}^N k(\cdot,x_n)\in H$ be the empirical mean embedding. Then 
$\Pr\left(\left\|\mu_{\hat{x}} - \mu_x\right\|_{H}  \le \frac{\sqrt{2B_k}}{\sqrt{N}}+\epsilon \right)\ge  1-e^{-\frac{\epsilon^2 N}{2B_k}}$,
or 
\begin{align*}
    \left\|\mu_{\hat{x}_i} - \mu_{x_i}\right\|_H &\le \frac{\sqrt{2B_k}}{\sqrt{N}} + \frac{\sqrt{2\alpha B_k}}{\sqrt{N}} = \frac{(1+\sqrt{\alpha})\sqrt{2B_k}}{\sqrt{N}} 
\end{align*}
with probability at least $1-e^{-\alpha}$, where $\alpha=\frac{\epsilon^2 N}{2B_k}$.

The proof will make use of the McDiarmid's inequality.

\begin{lemma}[\tb{McDiarmid's inequality} \cite{taylor04kernel}] \label{lemma:McDiarmid}
Let $x_1,\ldots,x_N\in\X$ be independent random variables and function $g\in\X^n\rightarrow \R$ be such that
$\sup_{u_1,\ldots,u_N,u_j'\in\X}\left|g(u_1,\ldots,u_N) - g(u_1,\ldots,u_{j-1},u_j',u_{j+1},\ldots,u_N)\right| \le c_j$
$\forall j=1,\ldots,N$. Then for all $\epsilon>0$ $\Pr\left( g(x_1,\ldots,x_N) - \E\left[g(x_1,\ldots,x_N)\right]\ge \epsilon\right) \le e^{-\frac{2\epsilon^2}{\sum_{n=1}^Nc_n^2}}$.
\end{lemma}

Namely, let $\phi(u)=k(\cdot,u)$, and thus $k(u,u)=\left\|\phi(u)\right\|_H^2$. Let us define
\begin{align*}
    g(S) &= \left\| \mu_{\hat{x}} -  \mu_x\right\|_H = \left\| \frac{1}{N}\sum_{n=1}^N \phi(x_n) - \mu_x \right\|_H,
\end{align*}
where $S=\{x_1,\ldots,x_N\}$ be the sample set. Define $S'=\{x_1,\ldots,x_{j-1},x_j',x_{j+1},\ldots,x_N\}$, i.e., let us replace in the sample set $x_j$ with $x_j'$. Then
\begin{align*}
    \left|g\left(S\right) - g\left(S'\right)\right|
    &= \left| \left\| \frac{1}{N}\sum_{n=1}^N \phi(x_n) - \mu_x \right\|_H
     - \left\| \frac{1}{N}\sum_{n=1; n\ne j}^N \phi(x_n) + \frac{1}{N}\phi(x_j') - \mu_x \right\|_H \right|\\
    &\le   \frac{1}{N} \left\| \phi(x_j) - \phi(x_j')\right\|_H \le \frac{1}{N} \left( \left\| \phi(x_j)\right\|_H + \left\| \phi(x_j')\right\|_H \right)
    \le \frac{1}{N} \left[ \sqrt{k(x_j,x_j)} + \sqrt{k\left(x_j',x_j'\right)} \right] \le \frac{2\sqrt{B_k}}{N}
\end{align*}
based on (i) the  reverse and the standard triangle inequality, and (ii) the boundedness of kernel $k$. By using the McDiarmid's inequality (Lemma~\ref{lemma:McDiarmid}), we get
\begin{align*}
 \Pr\left(g(S) - \E[g(S)] \ge \epsilon \right) &\le e^{-\frac{2\epsilon^2}{\sum_{n=1}^N\left( \frac{2\sqrt{B_k}}{N} \right)^2}}
    = e^{-\frac{2\epsilon^2}{N\frac{4B_k}{N^2}}} = e^{-\frac{\epsilon^2N}{2B_k}},
\end{align*}
or, in other words
\begin{align*}
 1 - e^{-\frac{\epsilon^2N}{2B_k}} \le \Pr\left(g(S)  < \E[g(S)] + \epsilon \right) \le \Pr\left(g(S)  \le \E[g(S)] + \epsilon \right).
\end{align*}
Considering the $\E[g(S)]$ term: since for a non-negative random variable ($a$) the 
$\E(a) = \E(a 1)\le \sqrt{\E (a^2)} \sqrt{\E (1^2)} = \sqrt{\E (a^2)}$ inequality holds due to the CBS, we obtain
\begin{align*}
	  \E[g(S)] & = \E\left[ \left\| \frac{1}{N}\sum_{n=1}^N \phi(x_n) - \mu_x \right\|_H \right]
	    \le \sqrt{\E\left[ \left\| \frac{1}{N}\sum_{n=1}^N \phi(x_n) - \mu_x \right\|_H^2 \right]}\\
	    &= \sqrt{\E\left[ \left<\frac{1}{N}\sum_{i=1}^N \phi(x_i) - \mu_x,\frac{1}{N}\sum_{j=1}^N \phi(x_j) - \mu_x \right>_H \right]}
	    = \sqrt{b+c+d}
\end{align*}
using that $\left\| a \right\|_H^2 = \sqrt{\left<a,a \right>_H}$. Here,
\begin{align*}
  b &= \E\left[\frac{1}{N^2}\left(\sum_{i,j=1; i\ne j}^Nk(x_i,x_j) + \sum_{i=1}^N k(x_i,x_i)\right)\right]
    = \frac{N(N-1)}{N^2}  \E_{t\sim x, t'\sim x}k(t,t') + \frac{N}{N^2} \E_{t\sim x} \left[k(t,t)\right],\\
  c &= - \frac{2}{N} \E\left[ \left<\sum_{i=1}^N \phi(x_i),\mu_x\right>_H\right] = - \frac{2N}{N} \E_{t\sim x, t'\sim x}\left[k(t,t')\right],\\
  d &= \E\left[\left\|\mu_x\right\|_H^2\right] = \E_{t\sim x, t'\sim x}\left[k(t,t')\right]
\end{align*}
applying the bilinearity of $\left<\cdot,\cdot \right>_H$, and the representation property of $\mu_x$.
Thus,
 \begin{align*}
  \sqrt{b+c+d}
  &= \sqrt{ \left[\frac{N-1}{N}-2+1\right] \E_{t\sim x, t'\sim x}\left[k(t,t')\right] +  \frac{1}{N}\E_{t\sim x} \left[k(t,t)\right] }\\
  &= \sqrt{\frac{1}{N} \left( \E_{t\sim x} \left[k(t,t)\right]  - \E_{t\sim x, t'\sim x}\left[k(t,t')\right]\right)}
  = \frac{\sqrt{\E_{t\sim x} \left[k(t,t)\right]  - \E_{t\sim x, t'\sim x}\left[k(t,t')\right]}}{\sqrt{N}}.
\end{align*}

Since 
  \begin{align*}
      \sqrt{\E_{t\sim x}\left[k(t,t)\right] - \E_{t\sim x,t'\sim x}\left[k(t,t')\right]}      
      &\le \sqrt{\left|\E_{t\sim x}\left[k(t,t)\right]\right|+\left|\E_{t\sim x,t'\sim x}\left[k(t,t')\right]\right|}
      \le \sqrt{\E_{t\sim x}\left|k(t,t)\right|+\E_{t\sim x,t'\sim x}\left|k(t,t')\right|},
  \end{align*}
  where we applied the triangle inequality, $\left|k(t,t)\right| = k(t,t)\le B_k$ and $\left|k(t,t')\right|\le \sqrt{k(t,t)} \sqrt{k(t',t')}$ (which holds to the CBS), we get
  $\sqrt{\E_{t\sim x}\left[k(t,t)\right] - \E_{t\sim x,t'\sim x}\left[k(t,t')\right]}  \le  \sqrt{B_k + \sqrt{B_k} \sqrt{B_k}}	= \sqrt{2B_k}$.
			  
To sum up, we obtained that $\left\|\mu_{x} - \mu_{\hat{x}}\right\|_H \le \frac{\sqrt{2B_k}}{\sqrt{N}} + \epsilon$ 
holds with probability at least $1-e^{-\frac{\epsilon^2 N}{2B_k}}$. This is what we wanted to prove.

\subsubsection{Proof of the Bound on $\| g_{\hat{\b{z}}} - g_{\b{z}}\|_{\H}^2$, $\| T_{\b{x}}-T_{\hat{\b{x}}}\|_{\L(\H)}^2$, $\|\sqrt{T} (T_{\hat{\b{x}}}+\lambda)^{-1}\|_{\L(\H)}^2$, $\|f_{\b{z}}^{\lambda}\|_{\H}^2$} 
Below, we present the detailed derivations of the upper bounds on $\| g_{\hat{\b{z}}} - g_{\b{z}}\|_{\H}^2$, $\| T_{\b{x}}-T_{\hat{\b{x}}}\|_{\L(\H)}^2$, $\|\sqrt{T} (T_{\hat{\b{x}}}+\lambda)^{-1}\|_{\L(\H)}^2$ and $\|f_{\b{z}}^{\lambda}\|_{\H}^2$.

\begin{compactitem}
      \item \tb{Bound on $\|  g_{\hat{\b{z}}} - g_{\b{z}} \|_{\H}^2$}: 	    By \eqref{eq:gz,gzhat}, we have $g_{\hat{\b{z}}} - g_{\b{z}} = \frac{1}{l}\sum_{i=1}^l \left[ K(\cdot,\mu_{\hat{x}_i}) - K(\cdot,\mu_{x_i}) \right] y_i$.
	    Applying Eq.~\eqref{eq:norm-eq}, the H{\"o}lder property of $K$, the homogenity of norms 
	    $\left\|a v\right\|=|a|\left\|v\right\|$ $(a\in\R)$,
	  	  assuming that  $y_i$ is bounded ($|y_i|\le C$),  and using \eqref{eq:emp-mean-emb-conv-rate}, we obtain
	    \begin{align*}
		    \left\|  g_{\hat{\b{z}}} - g_{\b{z}}  \right\|_{\H}^2  &\le  \frac{1}{l^2} l \sum_{i=1}^l \left\|  K(\cdot,\mu_{\hat{x}_i}) - K(\cdot,\mu_{x_i}) y_i \right\|_{\H}^2
		      \le \frac{L^2}{l} \sum_{i=1}^l y_i^2 \left\|\mu_{\hat{x}_i}-\mu_{x_i}\right\|_{H}^{2h}
		      \le \frac{L^2C^2}{l}\sum_{i=1}^l\left[ \frac{\left(1+\sqrt{\alpha}\right)\sqrt{2B_k}}{\sqrt{N}} \right]^{2h}\\
		      &=L^2 C^2  \frac{\left(1+\sqrt{\alpha}\right)^{2h} (2B_k)^h}{N^h} 
	    \end{align*}
	    with probability at least $1 - l e^{-\alpha}$, based on a union bound.

      \item \tb{Bound on $\| T_{\b{x}}-T_{\hat{\b{x}}}\|_{\L(\H)}^2$}: 	  Using the definition of $T_{\b{x}}$ and $T_{\hat{\b{x}}}$, and \eqref{eq:norm-eq} with the $\|\cdot\|_{\L(\H)}$ operator norm, we get
	  \begin{align}
	      \left\|T_{\b{x}} - T_{\hat{\b{x}}}\right\|_{\L(\H)}^2 &\le \frac{1}{l^2} l \sum_{i=1}^l \left\| T_{\mu_{x_i}} - T_{\mu_{\hat{x}_i}}\right\|_{\L(\H)}^2. \label{eq:Tx-Txi}
	  \end{align}
	  To upper bound the quantities $\| T_{\mu_{x_i}} - T_{\mu_{\hat{x}_i}}\|_{\L(\H)}^2$, let us see how $T_{\mu_u}$ acts
			\begin{align}
			    T_{\mu_u}(f) &=  K(\cdot,\mu_u)\delta_{\mu_u} (f) = K(\cdot,\mu_u) f(\mu_u).  \label{eq:T_u:act}
			\end{align}
	  If we can prove that
		\begin{align}
			\left\|(T_{\mu_u}-T_{\mu_v})(f)\right\|_{\H} &\le E \left\|f\right\|_{\H}, \label{eq:op-norm-estimation}
		\end{align}
		then this implies $\left\|T_{\mu_u}-T_{\mu_v}\right\|_{\L(\H)} \le E$. We continue with the l.h.s.\ of \eqref{eq:op-norm-estimation} using 
		  $\eqref{eq:T_u:act}$, \eqref{eq:norm-eq} with $n=2$, the homogenity of norms, the reproducing and H{\"o}lder property of $K$:
		  \begin{align*}
			\left\|(T_{\mu_u}-T_{\mu_v})(f)\right\|_{\H}^2
			  &=\left\| K(\cdot,\mu_u) \delta_{\mu_u}(f) - K(\cdot,\mu_v) \delta_{\mu_v}(f)\right\|_{\H}^2\\
			  &= \left\| K(\cdot,\mu_u) \left[\delta_{\mu_u}(f) -  \delta_{\mu_v}(f)\right]
			  +\left[ K(\cdot,\mu_u) - K(\cdot,\mu_v)\right] \delta_{\mu_v}(f)\right\|_{\H}^2 \nonumber\\
			  &\le  2 \left[ \left\| K(\cdot,\mu_u) \left[\delta_{\mu_u}(f) - \delta_{\mu_v}(f)\right]\right\|_{\H}^2
			  +\left\| K(\cdot,\mu_u) - K(\cdot,\mu_v) \delta_{\mu_v}(f) \right\|_{\H}^2\right]\\
			  &=  2 \left[ \left[\delta_{\mu_u}(f) - \delta_{\mu_v}(f)\right]^2 \left\| K(\cdot,\mu_u)\right\|_{\H}^2+
			  \left[ \delta_{\mu_v}(f) \right]^2 \left\| \left[K(\cdot,\mu_u)  - K(\cdot,\mu_v)\right] \right\|_{\H}^2\right]\\
			  &\le  2 \left[ \left[\delta_{\mu_u}(f) - \delta_{\mu_v}(f)\right]^2 K(\mu_u,\mu_u) + 
			  L^2\left[ \delta_{\mu_v}(f) \right]^2 \left\|\mu_u  - \mu_v \right\|_{H}^{2h}\right].
		  \end{align*}
		  By rewriting the first terms, we arrive at
		\begin{align*}
		      \delta_{\mu_u}(f) - \delta_{\mu_v}(f) &= \left<f,K(\cdot,\mu_u)\right>_{\H} - \left<f,K(\cdot,\mu_v)\right>_{\H}
		      \le \left|\left<f,K(\cdot,\mu_u)-K(\cdot,\mu_v)\right>_{\H}\right|\\
		      &\le \left\|f\right\|_{\H} \left\|K(\cdot,\mu_u)-K(\cdot,\mu_v)\right\|_{\H}  \le \left\|f\right\|_{\H}L\left\|\mu_u-\mu_v\right\|_H^h,\\
		      \delta_{\mu_v}(f) &= \left<f,K(\cdot,\mu_v)\right>_{\H} \le \left| \left<f,K(\cdot,\mu_v)\right>_{\H} \right|
		      \le \left\|f\right\|_{\H} \left\| K(\cdot,\mu_v) \right\|_{\H} = \left\|f\right\|_{\H} \sqrt{K(\mu_v,\mu_v)},
		\end{align*}
		where we applied the reproducing and H{\"o}lder property of $K$, the bilinearity of $\left<\cdot,\cdot\right>_{\H}$ and the CBS inequality. Hence
		\begin{align*}
		    \left\|(T_{\mu_u}-T_{\mu_v})(f)\right\|_{\H}^2 &\le 2\left[ \left\|f\right\|_{\H}^2 L^2\left\|\mu_u-\mu_v\right\|_H^{2h} K(\mu_u,\mu_u)
		    + L^2\left\|f\right\|_{\H}^2 K(\mu_v,\mu_v)\left\|\mu_u-\mu_v\right\|_H^{2h}\right]\\
		    &= 2 L^2\left\|f\right\|_{\H}^2 \left\|\mu_u-\mu_v\right\|_H^{2h}\left[K(\mu_u,\mu_u)+K(\mu_v,\mu_v)\right].
		\end{align*}
		Thus
	      \begin{align*}
		  E^2 = 2 L^2\|\mu_u-\mu_v\|_H^{2h}\left[K(\mu_u,\mu_u)+K(\mu_v,\mu_v)\right].
	      \end{align*}
		Exploiting this property in \eqref{eq:Tx-Txi}, \eqref{eq:bounded kernel}, and \eqref{eq:emp-mean-emb-conv-rate}
		\begin{align}
		      \left\|T_{\b{x}} - T_{\hat{\b{x}}}\right\|_{\L(\H)}^2
		      & \le  \frac{2L^2}{l}\sum_{i=1}^l\left\|\mu_{x_i}-\mu_{\hat{x}_i}\right\|_H^{2h}\left[K(\mu_{x_i},\mu_{x_i})+K(\mu_{\hat{x}_i},\mu_{\hat{x}_i})\right]
		      \le \frac{4B_K L^2}{l}\sum_{i=1}^l \frac{\left(1+\sqrt{\alpha}\right)^{2h} (2B_k)^h}{N^h}\nonumber\\
		      &= \frac{\left(1+\sqrt{\alpha}\right)^{2h} 2^{h+2}(B_k)^{h}B_K L^2}{N^h}.\label{eq:Tx-Txhat}
		\end{align}

      \item \tb{Bound on $\| \sqrt{T} (T_{\hat{\b{x}}}+\lambda)^{-1}\|_{\L(\H)}^2$}: First we rewrite $T_{\hat{\b{x}}}+\lambda$,
	      \begin{align*}
		  T_{\hat{\b{x}}}+\lambda &= (T + \lambda) - (T - T_{\hat{\b{x}}})
					  = \left[I- (T - T_{\hat{\b{x}}})(T+\lambda)^{-1}\right](T+\lambda).
	      \end{align*}
	      Let us now use the Neumann series of $I- (T - T_{\hat{\b{x}}})(T+\lambda)^{-1}$ 
	    \begin{eqnarray*}
		\sqrt{T} (T_{\hat{\b{x}}}+\lambda)^{-1} = \sqrt{T}(T+\lambda)^{-1} \sum_{n=0}^{\infty}\left[(T - T_{\hat{\b{x}}})(T+\lambda)^{-1}\right]^n 
	    \end{eqnarray*}
	      to have
	      \begin{align*}
		  \left\| \sqrt{T} (T_{\hat{\b{x}}}+\lambda)^{-1}\right\|_{\L(\H)} 
		    &=\left\| \sqrt{T}(T+\lambda)^{-1} \sum_{n=0}^{\infty}\left[(T - T_{\hat{\b{x}}})(T+\lambda)^{-1}\right]^n\right\|_{\L(\H)}\\
		    &\le\left\| \sqrt{T}(T+\lambda)^{-1}\right\|_{\L(\H)}  \left\|\sum_{n=0}^{\infty}\left[(T - T_{\hat{\b{x}}})(T+\lambda)^{-1}\right]^n\right\|_{\L(\H)}\\
		    &\le\left\| \sqrt{T}(T+\lambda)^{-1}\right\|_{\L(\H)}  \sum_{n=0}^{\infty}\left\|\left[(T - T_{\hat{\b{x}}})(T+\lambda)^{-1}\right]^n\right\|_{\L(\H)}\\
		    &\le \left\| \sqrt{T}(T+\lambda)^{-1}\right\|_{\L(\H)}  \sum_{n=0}^{\infty}\left\|(T - T_{\hat{\b{x}}})(T+\lambda)^{-1}\right\|_{\L(\H)}^n,    
	    \end{align*}
	    where $\left\|AB\right\|_{\L(\H)}\le \left\|A\right\|_{\L(\H)} \left\|B\right\|_{\L(\H)}$ and the triangle inequality was applied.
	    By the spectral theorem, the first term can be bounded as $\| \sqrt{T}(T+\lambda)^{-1}\|_{\L(\H)} \le \frac{1}{2\sqrt{\lambda}}$,	
	    whereas for the second term, applying a telescopic trick and a triangle inequality, we get
	    \begin{align*}
		\left\|(T - T_{\hat{\b{x}}})(T+\lambda)^{-1}\right\|_{\L(\H)}
		&= \left\|\left[(T - T_{\b{x}}) + (T_{\b{x}} - T_{\hat{\b{x}}})\right](T+\lambda)^{-1}\right\|_{\L(\H)}\\
		&\le \left\|(T - T_{\b{x}}) (T+\lambda)^{-1}\right\|_{\L(\H)} + \left\|(T_{\b{x}} - T_{\hat{\b{x}}})(T+\lambda)^{-1}\right\|_{\L(\H)}.
	    \end{align*}
	      We know that 
	      \begin{align}
		  \bm{\Theta}(\lambda,\b{z})&:=\|(T - T_{\b{x}}) (T+\lambda)^{-1}\|_{\L(\H)} \le \frac{1}{2} \label{eq:theta}
	      \end{align}
	      with probability at least $1-\frac{\eta}{3}$ \cite{caponnetto07optimal}.
	      Considering the second term, using \eqref{eq:Tx-Txhat} and $\|(T+\lambda)^{-1}\|_{\L(\H)}\le \frac{1}{\lambda}$ (by the spectral theorem),
		    \begin{align*}
			\left\|(T_{\b{x}} - T_{\hat{\b{x}}})(T+\lambda)^{-1}\right\|_{\L(\H)}
			& \le \left\|T_{\b{x}} - T_{\hat{\b{x}}}\right\|_{\L(\H)} \left\|(T+\lambda)^{-1}\right\|_{\L(\H)}
			\le \frac{\left(1+\sqrt{\alpha}\right)^{h} 2^{\frac{h}{2}+1}(B_k)^{\frac{h}{2}}(B_K)^{\frac{1}{2}} L}{N^{\frac{h}{2}}}  \frac{1}{\lambda}.
		    \end{align*}
		    For fixed $\lambda$, the value of $N$ can be chosen such that
		    \begin{align}
			  \frac{\left(1+\sqrt{\alpha}\right)^{h} 2^{\frac{h}{2}+1}(B_k)^{\frac{h}{2}}(B_K)^{\frac{1}{2}} L}{N^{\frac{h}{2}}}  \frac{1}{\lambda} &\le \frac{1}{4}\Leftrightarrow
			  \frac{\left(1+\sqrt{\alpha}\right)^{h} 2^{\frac{h}{2}+3}(B_k)^{\frac{h}{2}}(B_K)^{\frac{1}{2}} L}{\lambda}  \le N^{\frac{h}{2}}\Leftrightarrow \nonumber\\
			  \frac{\left(1+\sqrt{\alpha}\right)^{2} 2^{\frac{h+6}{h}}B_k (B_K)^{\frac{1}{h}} L^{\frac{2}{h}}}{\lambda^{\frac{2}{h}}}  &\le N \label{eq:N}.
		    \end{align}
	    In this case $\left\|(T - T_{\hat{\b{x}}})(T+\lambda)^{-1}\right\|_{\L(\H)}\le \frac{3}{4}$ (the Neumann series trick is legitimate) and
	    \begin{align}
		\left\| \sqrt{T} (T_{\hat{\b{x}}}+\lambda)^{-1}\right\|_{\L(\H)} &\le \frac{1}{2\sqrt{\lambda}} \frac{1}{1-\frac{3}{4}} \le \frac{2}{\sqrt{\lambda}}. \label{eq:common term in S-1 S0}
	    \end{align}
      \item \tb{Bound on $\|f_{\b{z}}^{\lambda}\|_{\H}^2$}: 	    Using the explicit form of $f_{\b{z}}^{\lambda}$ [\eqref{eq:f_zlambda}], \eqref{eq:opnorm-fit}, the positivity of $T_{\b{x}}$ $\left[\Rightarrow\left\|(T_{\b{x}}+\lambda)^{-1}\right\|_{\L(\H)}\le \frac{1}{\lambda}\right]$, the homogenity of norms, Eq.~\eqref{eq:norm-eq}, the boundedness assumption on $y_i$ ($|y_i| \le C$),
	    the reproducing property and the boundedness of $K$ [Eq.~\eqref{eq:bounded kernel}], we get
	    \begin{align*}
		  \left\|f_{\b{z}}^{\lambda}\right\|_{\H} \le \left\| (T_{\b{x}}+\lambda)^{-1}\right\|_{\L(\H)} \left\|g_{\b{z}} \right\|_{\H} \le  \frac{1}{\lambda} \left\|g_{\b{z}} \right\|_{\H},
	    \end{align*}
	    where
	    \begin{align*}
		  \left\|g_{\b{z}} \right\|_{\H}^2& \le \frac{1}{l^2} l \sum_{i=1}^l\left\| K(\cdot,\mu_{x_i}) y_i \right\|_{\H}^2
		  \le \frac{1}{l} \sum_{i=1}^l C^2 \left\|K(\cdot,\mu_{x_i})\right\|_{\H}^2
		  = \frac{1}{l} \sum_{i=1}^l C^2 K(\mu_{x_i},\mu_{x_i}) \le \frac{1}{l} \sum_{i=1}^l C^2 B_K = C^2 B_K.
	    \end{align*}
	    Thus, we have obtained that
	    \begin{align}
		\left\|f_{\b{z}}^{\lambda}\right\|_{\H}^2& \le \frac{1}{\lambda^2}C^2 B_K. \label{eq:bound on fzlambda}
	    \end{align}

\end{compactitem}

\subsubsection{Final Step of the Proof (Union Bound)}
    Until now, we obtained that if 
    \begin{compactenum}
	  \item the sample number $N$ satisfies Eq.~\eqref{eq:N},  
	  \item \eqref{eq:emp-mean-emb-conv-rate} holds for $\forall i=1,\ldots,l$ (which has probability at least $1 - l e^{-\alpha} = 1 - e^{-[\alpha-\log(l)]}= 1-e^{-\delta}$ applying a union bound argument; $\alpha=\log(l)+\delta$), and
	  \item $\bm{\Theta}(\lambda,\b{z})\le \frac{1}{2}$ is fulfilled [see Eq.~\eqref{eq:theta}], then
    \end{compactenum}
  \begin{align*}
	S_{-1} + S_{0} & \le \frac{4}{\lambda} \left[L^2 C^2  \frac{\left(1+\sqrt{\alpha}\right)^{2h} (2B_k)^h}{N^h}
	       + \frac{\left(1+\sqrt{\alpha}\right)^{2h} 2^{h+2}(B_k)^{h}B_K L^2}{N^h}  \frac{C^2 B_K}{\lambda^2}\right]\\
	      & = \frac{4L^2 C^2 \left(1+\sqrt{\alpha}\right)^{2h} (2B_k)^h}{\lambda N^h} \left[  1 + \frac{4(B_K)^2}{\lambda^2}\right].
  \end{align*}
  By taking into account \cite{caponnetto07optimal}'s bounds for $S_1$ and $S_2$
  \begin{align*}
    S_1 &\le 32 \log^2\left(\frac{6}{\eta}\right)\left[\frac{B_K M^2}{l^2\lambda} + \frac{\Sigma^2 \N(\lambda)}{l} \right], &
    S_2 &\le 8 \log^2\left(\frac{6}{\eta}\right) \left[\frac{4B_K^2\B(\lambda)}{l^2\lambda} + \frac{B_K \A(\lambda)}{l\lambda} \right],
  \end{align*}
  plugging all the expressions to \eqref{eq:5-bound}, we obtain Theorem~\ref{theo1} via a union bound.

\subsubsection{Proof of Consequence~\ref{conseq:conv-rate}}
Since constant multipliers do not matter in the orders of rates, we discard them in the (in)equalities below.
Our goal is to choose $\lambda=\lambda_{l,N}$ such that 
\begin{compactitem}
  \item $\lim_{l,N\rightarrow\infty}\lambda_{l,N}=0$, and
  \item in Theorem~\ref{theo1}: (i') $\frac{\log(l)}{\lambda^{\frac{2}{h}}} \le N$, (i) $l \lambda^{\frac{b+1}{b}}\ge 1$,\footnote{ 
$\N(\lambda)$ can be upper bounded by (constant multipliers are discarded) $\lambda^{-\frac{1}{b}}$ \cite{caponnetto07optimal}. Using this upper bound in the $l$ constraint of Theorem~\ref{theo1} we get $l\ge\frac{\lambda^{-\frac{1}{b}}}{\lambda}\Leftrightarrow l \lambda^{\frac{1}{b}+1=\frac{b+1}{b}}\ge 1$. } and (ii) $r(l,N,\lambda) = \frac{\log^h(l)}{N^h\lambda^3}+\lambda^c + \frac{\lambda^{c-2}}{l^2} + \frac{\lambda^{c-1}}{l} + \frac{1}{l^2\lambda} + \frac{1}{l\lambda^{\frac{1}{b}}}\rightarrow 0$.
\end{compactitem}
In $r(l,N,\lambda)$ we will require that the first term goes to zero $\left[\frac{\log^h(l)}{N^h\lambda^3}\rightarrow 0\right]$, which implies 
$\frac{\log(l)}{N\lambda^{\frac{3}{h}}}\rightarrow 0$ and hence $\frac{\log(l)}{N\lambda^{\frac{2}{h}}}\rightarrow 0$. Thus constraint (i') can be discarded, 
and our goal is to fulfill (i)-(ii). Since
\begin{compactenum}
	    \item $2-c\le 1$ ($\Leftrightarrow$ $1\le c$), $\frac{1}{l^2\lambda^{2-c}} = \frac{\lambda^{c-2}}{l^2} \le \frac{1}{l^2\lambda}$ (in order), and
	    \item $c-1 \ge 0$ ($\Leftrightarrow$ $1\le c$), $\frac{\lambda^{c-1}}{l} \le \frac{1}{l\lambda^{\frac{1}{b}}}$ (in order)
\end{compactenum}
	  condition (i)-(ii) reduces to
	  \begin{align}
	      r(l,N,\lambda) &= \frac{\log^h(l)}{N^h\lambda^3}+\lambda^c + \frac{1}{l^2\lambda} + \frac{1}{l\lambda^{\frac{1}{b}}}\rightarrow 0, \text{ subject to } l \lambda^{\frac{b+1}{b}}\ge 1\label{eq:f1}.
      \end{align}
Our goal is to study the behavior of this quantity in terms of the $(l,N,\lambda)$ triplet; $1<b$, $c\in [1,2]$, $h\in(0,1]$. To do so, we  
\begin{compactenum}
  \item choose $\lambda$ such a way that two terms match in order (and $\lambda=\lambda_{l,N}\rightarrow 0$);
  \item setting $l=N^a$ ($a>0$) we examine under what conditions (i)-(ii) the convergence of $r$ to $0$ holds with the constraint $l \lambda^{\frac{b+1}{b}}\ge 1$ satisfied, (iii) are the matched terms also dominant, i.e., give the convergence rate.
\end{compactenum}
We carry out the computation for all the $\binom{4}{2}=6$ pairs in Eq.~\eqref{eq:f1}. Below we give the derivation of the results summarized in
Table~\ref{tab:conv-rates}.

\begin{compactitem}
  \item $\boxed{1} = \boxed{2}$ in Eq.~\eqref{eq:f1} [i.e., the first and second terms are equal in Eq.~\eqref{eq:f1}]:
      \begin{compactitem}
	  \item (i)-(ii): Exploiting $\frac{h}{c+3}>0$ in the $\lambda$ choice, we get
	      \begin{align}
		\frac{\log^h(l)}{N^h\lambda^3} &=\lambda^c \Leftrightarrow \left[\frac{\log(l)}{N}\right]^h = \lambda^{c+3} \Leftrightarrow \left[\frac{\log(l)}{N}\right]^{\frac{h}{c+3}}  = \lambda \rightarrow 0, \text{ if } \frac{\log(l)}{N}\rightarrow 0. \nonumber\\
		r(l,N) &= \left[\frac{\log(l)}{N}\right]^{\frac{hc}{c+3}} + \frac{1}{l^2 \left[\frac{\log(l)}{N}\right]^{\frac{h}{c+3}}} + \frac{1}{l\left[\frac{\log(l)}{N}\right]^{\frac{h}{b(c+3)}}}. \nonumber\\
		  r(N) &= \left[ \frac{\log(N)}{N}\right]^{\frac{hc}{c+3}} + \frac{1}{N^{2a} \left[\frac{\log(N)}{N}\right]^{\frac{h}{c+3}}} + \frac{1}{N^a \left[\frac{\log(N)}{N}\right]^{\frac{h}{b(c+3)}}}\nonumber\\
		       &= \left[ \frac{\log(N)}{N}\right]^{\frac{hc}{c+3}} + \frac{N^{\frac{h}{c+3}}}{N^{2a} \log^{\frac{h}{c+3}}(N)} + \frac{N^{\frac{h}{b(c+3)}}}{N^a \log^{\frac{h}{b(c+3)}}(N)}. \label{eq:fN1}
	      \end{align}
	      Here, 
		\begin{compactitem}
		    \item (ii): $r(N)\rightarrow 0$ if 
			\begin{compactitem}
			    \item $\ovalbox{1}\rightarrow 0$: [i.e., the first term goes to zero in Eq.~\eqref{eq:fN1}]; no constraint using that $\frac{hc}{c+3}>0$.
			    \item $\ovalbox{2}\rightarrow 0$: $2a\ge \frac{h}{c+3}$ [$\Leftarrow \frac{h}{c+3}>0$].
			    \item $\ovalbox{3}\rightarrow 0$: $a\ge \frac{h}{b(c+3)}$ [$\Leftarrow \frac{h}{b(c+3)}>0$],
			\end{compactitem}
			i.e., $a\ge \max\left(\frac{h}{2(c+3)},\frac{h}{b(c+3)}\right) = \frac{h}{(c+3)\min(2,b)}$.
		    \item (i): We require $N^a\left(\left[\frac{\log(N)}{N}\right]^{\frac{h}{c+3}}\right)^{\frac{b+1}{b}}\ge 1 \Leftrightarrow \frac{\log^{\frac{h}{c+3}\frac{b+1}{b}}(N)}{N^{\frac{h}{c+3}\frac{b+1}{b}-a}} \ge 1$.
			  Since $\frac{h}{c+3}\frac{b+1}{b}>0$, it is sufficient to have $\frac{h}{c+3}\frac{b+1}{b}-a\le 0 \Leftrightarrow \frac{h(b+1)}{(c+3)b}\le a$.
		\end{compactitem}
		To sum up, for (i)-(ii) we got $a\ge \max\left(\frac{h}{(c+3)\min(2,b)},\frac{h(b+1)}{(c+3)b}\right)$.
	  \item (iii): 
	    \begin{compactitem}
		\item (i): $\frac{h(b+1)}{(c+3)b}\le a$.
		\item $\ovalbox{1}\rightarrow 0$: no constraint.
		\item $\ovalbox{1} \ge \ovalbox{2}$ [i.e., the first term dominates the second one in Eq.~\eqref{eq:fN1}]:
			$\left[ \frac{\log(N)}{N}\right]^{\frac{hc}{c+3}} \ge  \frac{N^{\frac{h}{c+3}}}{N^{2a} \log^{\frac{h}{c+3}}(N)} \Leftrightarrow \log^{\frac{hc}{c+3}+\frac{h}{c+3}}(N) \ge N^{\frac{hc}{c+3}+\frac{h}{c+3}-2a}$.
		    Thus, since $\frac{h(c+1)}{c+3}>0$ we need $\frac{h(c+1)}{c+3}-2a \le 0$, i.e., $\frac{h(c+1)}{2(c+3)} \le a$.
		\item $\ovalbox{1} \ge \ovalbox{3}$ [i.e., the first term dominates the third one in Eq.~\eqref{eq:fN1}]: $\left[ \frac{\log(N)}{N}\right]^{\frac{hc}{c+3}} \ge  \frac{N^{\frac{h}{b(c+3)}}}{N^a \log^{\frac{h}{b(c+3)}}(N)} \Leftrightarrow  \log^{\frac{hc}{c+3}+\frac{h}{b(c+3)}}(N) \ge N^{\frac{h}{b(c+3)} + \frac{hc}{c+3}-a}$. 
		    Since $\frac{hc}{c+3}+\frac{h}{b(c+3)}>0$ we require $\frac{h}{b(c+3)} + \frac{hc}{c+3}-a\le 0$, i.e., $\frac{h}{b(c+3)} + \frac{hc}{c+3} \le a$.
	    \end{compactitem}
	    To sum up, the obtained condition for $a$ is $\max\left( \frac{h}{b(c+3)} + \frac{hc}{c+3}, \frac{h(c+1)}{2(c+3)}\right) = \frac{h\max\left( \frac{1}{b} + c, \frac{c+1}{2}\right)}{c+3}  \le a$.
	    Since $\frac{1}{b} + c \ge \frac{c+1}{2}  \Leftrightarrow \frac{1}{b} + \frac{c}{2} \ge \frac{1}{2} [\Leftarrow c\ge 1, b>0]$,
	    we got
	    \begin{align*}
		\max\left(\frac{h\left(\frac{1}{b}+c\right)}{c+3},\frac{h(b+1)}{(c+3)b}\right)&\le a, & r(N)&= \left[ \frac{\log(N)}{N}\right]^{\frac{hc}{c+3}} \rightarrow 0.
	    \end{align*}
      \end{compactitem}
  \item $\boxed{1} = \boxed{3}$ in Eq.~\eqref{eq:f1}:
      \begin{compactitem}
	  \item (i)-(ii): Using in the $\lambda$ choice that $\frac{h}{2}>0$, we obtain that
		    \begin{align*}
			 \frac{\log^h(l)}{N^h\lambda^3} &= \frac{1}{l^2\lambda} \Leftrightarrow \frac{l^2\log^h(l)}{N^h} = \lambda^2 \Leftrightarrow \frac{l\log^{\frac{h}{2}}(l)}{N^{\frac{h}{2}}} = \lambda \rightarrow 0\text{, if } a<\frac{h}{2} \text{ in } l=N^a.\\
			  r(l,N) &= \left[\frac{l\log^{\frac{h}{2}}(l)}{N^{\frac{h}{2}}}\right]^c + \frac{1}{l^2 \frac{l\log^{\frac{h}{2}}(l)}{N^{\frac{h}{2}}}} + \frac{1}{l \left[\frac{l\log^{\frac{h}{2}}(l)}{N^{\frac{h}{2}}}\right]^{\frac{1}{b}}}.\\
			  r(N)   &= N^{ac-\frac{hc}{2}} \log^{\frac{hc}{2}}(N) + \frac{1}{N^{3a-\frac{h}{2}} \log^{\frac{h}{2}}(N)} + \frac{1}{N^{a+\frac{a}{b}-\frac{h}{2b}} \log^{\frac{h}{2b}}(N)}.
		    \end{align*}
		Here, 
	      \begin{compactitem}
		  \item (ii): $r(N)\rightarrow 0$ if 
			\begin{compactitem}
			      \item $\ovalbox{1}\rightarrow 0$: $ac-\frac{hc}{2}=c\left(a-\frac{h}{2}\right)<0$ [$\Leftarrow$ $\frac{hc}{2}>0$], i.e., $a<\frac{h}{2}$ using that $c>0$. 
			      \item $\ovalbox{2}\rightarrow 0$: $3a-\frac{h}{2}\ge 0$ [$\Leftarrow$ $\frac{h}{2}>0$], i.e., $\frac{h}{6}\le a$.
			      \item $\ovalbox{3}\rightarrow 0$: $a+\frac{a}{b}-\frac{h}{2b} \ge 0$ [$\Leftarrow$ $\frac{h}{2b}>0$], i.e., $\frac{h}{2b\left(1+\frac{1}{b}\right)} = \frac{h}{2b\frac{b+1}{b}}=\frac{h}{2(b+1)}\le a$ exploiting that $1+\frac{1}{b}>0$.
			\end{compactitem}
		      In other words,  the requirement is $\max\left( \frac{h}{6},\frac{h}{2(b+1)}\right) \le a < \frac{h}{2}$.
		  \item (i): $N^a \left[\frac{N^a\log^{\frac{h}{2}}(N)}{N^{\frac{h}{2}}}\right]^{\frac{b+1}{b}}\ge 1 \Leftrightarrow \frac{\log^{\frac{h}{2}\frac{b+1}{b}}(N)}{N^{\frac{h}{2}\frac{b+1}{b}-a-a\frac{b+1}{b}}} \ge 1$.
		      Since $\frac{h}{2}\frac{b+1}{b}>0$ it is enough to have $\frac{h}{2}\frac{b+1}{b}-a-a\frac{b+1}{b} \le 0  \Leftrightarrow \frac{h}{2}\frac{b+1}{b} \le a\left(1+\frac{b+1}{b}\right) = a \frac{2b+1}{b} \Leftrightarrow \frac{h}{2}\frac{b+1}{2b+1}\le a$
		      using that $2b+1>0$, $b>0$ [$\Leftarrow b>1$].
	      \end{compactitem}
	      To sum up, for (i)-(ii) we obtained $\max\left( \frac{h}{6},\frac{h}{2(b+1)},\frac{h}{2}\frac{b+1}{2b+1}\right) \le a < \frac{h}{2}$.
	  \item (iii): 
	    \begin{compactitem}
	      \item (i): $\frac{h}{2}\frac{b+1}{2b+1}\le a$
	      \item $\ovalbox{2}\rightarrow 0$: $\frac{h}{6}\le a$.
		\item $\ovalbox{2} \ge \ovalbox{1}$: $\frac{1}{N^{3a-\frac{h}{2}} \log^{\frac{h}{2}}(N)} \ge N^{ac-\frac{hc}{2}} \log^{\frac{hc}{2}}(N) \Leftrightarrow N^{\frac{h}{2}-3a+\frac{hc}{2}-ac}\ge \log^{\frac{h(c+1)}{2}}(N)$.
		    Thus, since $\frac{h(c+1)}{2}>0$ we need $\frac{h}{2}-3a+\frac{hc}{2}-ac>0$, i.e., $\frac{h(c+1)}{2(c+3)}=\frac{h(c+3-2)}{2(c+3)}=\frac{h}{2}-\frac{h}{c+3}>a$, using that $c+3>0$.
		\item $\ovalbox{2} \ge \ovalbox{3}$:
			$\frac{1}{N^{3a-\frac{h}{2}} \log^{\frac{h}{2}}(N)} \ge \frac{1}{N^{a+\frac{a}{b}-\frac{h}{2b}} \log^{\frac{h}{2b}}(N)} \Leftrightarrow N^{a+\frac{a}{b}-\frac{h}{2b}+\frac{h}{2}-3a} \ge \log^{\frac{h}{2}-\frac{h}{2b}}(N)$.
		    Since $\frac{h}{2}-\frac{h}{2b} = \frac{h}{2}\left(1-\frac{1}{b}\right)>0$ using that $h>0$ and $b>1$, we need $a+\frac{a}{b}-\frac{h}{2b}+\frac{h}{2}-3a>0$, i.e., 
                    $a\left(1+\frac{1}{b}-3\right)>\frac{h}{2}\left(\frac{1}{b}-1\right) \Leftrightarrow a \left(\frac{1}{b}-2\right) >\frac{h}{2}\left(\frac{1}{b}-1\right)$.
		  Using that $b>1$, $0>\frac{1}{b}-1>\frac{1}{b}-2$; hence $a < \frac{\frac{h}{2}\left(\frac{1}{b}-1\right)}{\frac{1}{b}-2}$.
	    \end{compactitem}
		  To sum up, we got
		  \begin{align*}
		      \max\left(\frac{h}{6},\frac{h}{2}\frac{b+1}{2b+1}\right)&\le a < \min\left( \frac{h}{2}-\frac{h}{c+3}, \frac{\frac{h}{2}\left(\frac{1}{b}-1\right)}{\frac{1}{b}-2} \right) &
		      r(N) &= \frac{1}{N^{3a-\frac{h}{2}} \log^{\frac{h}{2}}(N)} \rightarrow 0.
		  \end{align*}
      \end{compactitem}
  \item $\boxed{1} = \boxed{4}$ in Eq.~\eqref{eq:f1}:
      \begin{compactitem}
	  \item (i)-(ii): Using in the $\lambda$ choice that $\frac{b}{3b-1}>0$, we get
		\begin{align*}
		  \frac{\log^h(l)}{N^h\lambda^3} &= \frac{1}{l\lambda^{\frac{1}{b}}} \Leftrightarrow \frac{l\log^h(l)}{N^h} = \lambda^{3-\frac{1}{b}=\frac{3b-1}{b}} \Leftrightarrow \left[\frac{l\log^h(l)}{N^h}\right]^{\frac{b}{3b-1}} = \lambda\rightarrow 0\text{, if $h>a$ in $l=N^a$}.\\
		  r(l,N) &= \left[\frac{l\log^h(l)}{N^h}\right]^{\frac{bc}{3b-1}} + \frac{1}{l^2 \left[\frac{l\log^h(l)}{N^h}\right]^{\frac{b}{3b-1}}} + \frac{1}{l \left[\frac{l\log^h(l)}{N^h}\right]^{\frac{1}{3b-1}}}.\\
		  r(N) &= \left[\frac{\log^h(N)}{N^{h-a}}\right]^{\frac{bc}{3b-1}} + \frac{1}{N^{2a+\frac{ab}{3b-1}-\frac{hb}{3b-1}}\log^{\frac{hb}{3b-1}}(N)} + \frac{1}{N^{a+\frac{a}{3b-1}-\frac{h}{3b-1}}\log^{\frac{h}{3b-1}}(N)}.
		\end{align*}
		Here, 
		\begin{compactitem}
		    \item (ii):
			  $r(N)\rightarrow 0$, if 
			  \begin{compactitem}
			      \item $\ovalbox{1}\rightarrow 0$: $h-a>0$ using that $h>0$ and $\frac{bc}{3b-1}>0$, i.e., $a<h$,
			      \item $\ovalbox{2}\rightarrow 0$: $2a+\frac{ab}{3b-1}-\frac{hb}{3b-1}\ge 0$ [using that $\frac{hb}{3b-1}>0$]. In other words,
				  $a\left(2+\frac{b}{3b-1}\right)  \ge \frac{hb}{3b-1} \Leftrightarrow a \ge \frac{\frac{hb}{3b-1}}{\left(2+\frac{b}{3b-1}\right)}=\frac{hb}{3b-1}\frac{3b-1}{6b-2+b}=\frac{hb}{7b-2}$
				  using that $\left(2+\frac{b}{3b-1}\right)>0$.
			      \item $\ovalbox{3}\rightarrow 0$: $a+\frac{a}{3b-1}-\frac{h}{3b-1}\ge 0$ [using that $\frac{h}{3b-1}>0$], i.e., 
				    $a\left(1+\frac{1}{3b-1}\right) \ge \frac{h}{3b-1}\Leftrightarrow a \ge \frac{\frac{h}{3b-1}}{1+\frac{1}{3b-1}} = \frac{h}{3b-1}\frac{3b-1}{3b-1+1}=\frac{h}{3b}$
				  making use of $\left(1+\frac{1}{3b-1}\right)>0$.
			  \end{compactitem}
			  Thus, we need $\max\left(\frac{hb}{7b-2},\frac{h}{3b}\right)\le a<h$.
		    \item (i): $N^a \left(\left[\frac{N^a\log^h(N)}{N^h}\right]^{\frac{b}{3b-1}}\right)^{\frac{b+1}{b}} \ge 1 \Leftrightarrow \frac{\log^{\frac{h(b+1)}{3b-1}}(N)}{N^{\frac{h(b+1)}{3b-1}-a-a\frac{b+1}{3b-1}}} \le 1$.
		      Since $\frac{h(b+1)}{3b-1}>0$, it is sufficient $\frac{h(b+1)}{3b-1}-a-a\frac{b+1}{3b-1}\le 0 \Leftrightarrow \frac{h(b+1)}{3b-1} \le a \left(1+\frac{b+1}{3b-1}\right) = a \frac{3b-1+b+1}{3b-1} = a \frac{4b}{3b-1}
			      \Leftrightarrow \frac{h(b+1)}{4b} \le a$,
		      where we used that $4b>0$, $3b-1>0$ [$\Leftarrow b>1$].
		\end{compactitem}
		To sum up, for (i)-(ii) we received $\max\left(\frac{hb}{7b-2},\frac{h}{3b},\frac{h(b+1)}{4b}\right)\le a<h$.
	  \item (iii): 
	    \begin{compactitem}
		\item (i): $\frac{h(b+1)}{4b}\le a$.
		\item $\ovalbox{3}\rightarrow 0$: $a \ge \frac{h}{3b}$. 
		\item $\ovalbox{3} \ge \ovalbox{1}$: $\frac{1}{N^{a+\frac{a}{3b-1}-\frac{h}{3b-1}}\log^{\frac{h}{3b-1}}(N)} \ge \left[\frac{\log^h(N)}{N^{h-a}}\right]^{\frac{bc}{3b-1}} \Leftrightarrow N^{\frac{(h-a)bc}{3b-1}-a-\frac{a}{3b-1}+\frac{h}{3b-1}}\ge \log^{\frac{h(bc+1)}{3b-1}}(N)$.
		    Since $\frac{h(bc+1)}{3b-1}>0$, we need $\frac{(h-a)bc}{3b-1}-a-\frac{a}{3b-1}+\frac{h}{3b-1} > 0 \Leftrightarrow \frac{h(bc+1)}{3b-1} > a\left(\frac{bc}{3b-1}+1+\frac{1}{3b-1}\right)
			 \Leftrightarrow \frac{h(bc+1)}{3b-1} > a\left(1+\frac{bc+1}{3b-1}\right)
			\Leftrightarrow \frac{h(bc+1)}{3b-1} > a\frac{3b-1+bc+1}{3b-1}
			 \Leftrightarrow \frac{h(bc+1)}{3b-1} > a\frac{3b+bc}{3b-1}
			\Leftrightarrow \frac{h(bc+1)}{3b+bc}>a$

		    using at the last step that $3b-1>0$ and $3b+bc>0$.
		\item $\ovalbox{3} \ge \ovalbox{2}$: $\frac{1}{N^{a+\frac{a}{3b-1}-\frac{h}{3b-1}}\log^{\frac{h}{3b-1}}(N)} \ge \frac{1}{N^{2a+\frac{ab}{3b-1}-\frac{hb}{3b-1}}\log^{\frac{hb}{3b-1}}(N)} \Leftrightarrow\\
			\log^{\frac{h(b-1)}{3b-1}}(N)\ge N^{-2a-\frac{ab}{3b-1}+\frac{hb}{3b-1} +a +\frac{a}{3b-1}-\frac{h}{3b-1} }$. 
		    Since $\frac{h(b-1)}{3b-1}>0$, we require that 
		    $-2a-\frac{ab}{3b-1}+\frac{hb}{3b-1} +a +\frac{a}{3b-1}-\frac{h}{3b-1} \le 0  \Leftrightarrow \frac{h(b-1)}{3b-1} \le a \left(1+\frac{b-1}{3b-1}\right)
	    	     \Leftrightarrow \frac{h(b-1)}{3b-1} \le a \frac{3b-1+b-1}{3b-1}
		     \Leftrightarrow \frac{h(b-1)}{4b-2} \le a$ using that $3b-1>0$ and $4b-2>0$.
	    \end{compactitem}
      \end{compactitem}
      To sum up, we obtained that
      \begin{align*}
	    \max\left(\frac{h(b-1)}{4b-2},\frac{h}{3b},\frac{h(b+1)}{4b}\right) &\le a < \frac{h(bc+1)}{3b+bc}, &
	    r(N)&=\frac{1}{N^{a+\frac{a}{3b-1}-\frac{h}{3b-1}}\log^{\frac{h}{3b-1}}(N)} \rightarrow 0.
      \end{align*}
  \item $\boxed{2} = \boxed{3}$ in Eq.~\eqref{eq:f1}:
   \begin{compactitem}
	  \item (i)-(ii):
	      \begin{align*}
		\lambda^c &= \frac{1}{l^2\lambda} \Leftrightarrow \lambda^{c+1} = \frac{1}{l^2} \Leftrightarrow \lambda = \frac{1}{l^{\frac{2}{c+1}}} \rightarrow 0, \text{ if }l\rightarrow \infty. \quad[\Leftarrow \frac{2}{c+1}>0]\\
		r(l,N) &= \frac{l^{\frac{6}{c+1}}\log^h(l)}{N^h} + \frac{1}{l^{\frac{2c}{c+1}}} + \frac{l^{\frac{2}{b(c+1)}}}{l} \Rightarrow
		r(N) = \frac{\log^h(N)}{N^{h-\frac{6a}{c+1}}} + \frac{1}{N^{\frac{2ac}{c+1}}} + \frac{1}{N^{a\left(1-\frac{2}{b(c+1)}\right)}}.
	      \end{align*}
	      Here, 
	      \begin{compactitem}
		  \item (ii):
			$r(N)\rightarrow 0$ if 
			\begin{compactitem}
			    \item $\ovalbox{1}\rightarrow 0$: $h-\frac{6a}{c+1}>0$ since $h>0$, i.e., $a<\frac{h(c+1)}{6}$ using that $c+1>0$.
			    \item $\ovalbox{2}\rightarrow 0$: $\frac{2ac}{c+1}>0$ -- this condition is satisfied by our assumptions ($a>0$, $c>0$).
			    \item $\ovalbox{3}\rightarrow 0$: $a\left(1-\frac{2}{b(c+1)}\right)>0$. Using that $a>0$, $b>0$, $c+1>0$ this 
				  requirement is $1 >\frac{2}{b(c+1)}\Leftrightarrow b(c+1)>2 [\Leftarrow b>1, c\ge 1]$.
			\end{compactitem}
		      Thus, we need $a <\frac{h(c+1)}{6}$.
		  \item (i): $N^a \left(\frac{1}{N^{\frac{2a}{c+1}}}\right)^{\frac{b+1}{b}} \ge 1 \Leftrightarrow N^{a-\frac{2a(b+1)}{(c+1)b}} \ge 1$.
		      Thus it is enough to satisfy $a-\frac{2a(b+1)}{(c+1)b} > 0\Leftrightarrow 1 > \frac{2(b+1)}{(c+1)b}$,
		      where we used that $a>0$.
	      \end{compactitem}
	      To sum up, for (i)-(ii) we obtained $a < \frac{h(c+1)}{6}$, $1 > \frac{2(b+1)}{(c+1)b}$.
	  \item (iii): 
	    \begin{compactitem}
		\item (i): $1 > \frac{2(b+1)}{(c+1)b}$.
		\item $\ovalbox{2}\rightarrow 0$: no constraint.
		\item $\ovalbox{2} \ge \ovalbox{1}$: $\frac{1}{N^{\frac{2ac}{c+1}}} \ge \frac{\log^h(N)}{N^{h-\frac{6a}{c+1}}} \Leftrightarrow N^{h-\frac{6a}{c+1}-\frac{2ac}{c+1}} \ge \log^h(N)$.
		    Thus, since $h>0$ we require that $h-\frac{6a}{c+1}-\frac{2ac}{c+1} > 0 \Leftrightarrow h>a\frac{6+2c}{c+1} \Leftrightarrow \frac{h(c+1)}{6+2c} > a$,
		     where the $6+2c>0$, $c+1>0$ relations were exploited [$\Leftarrow c \ge 1$].
		\item $\ovalbox{2} \ge \ovalbox{3}$: $\frac{1}{N^{\frac{2ac}{c+1}}} \ge \frac{1}{N^{a\left(1-\frac{2}{b(c+1)}\right)}} \Leftrightarrow N^{ a\left(1-\frac{2}{b(c+1)}\right) - \frac{2ac}{c+1}} \ge 1$.
		    Hence, by $a>0$ and $c+1>0$ we need $a\left(1-\frac{2}{b(c+1)}\right) - \frac{2ac}{c+1} >0 \Leftrightarrow
			a \frac{b(c+1)-2}{b(c+1)} > \frac{2ac}{c+1}\Leftrightarrow
			b(c+1)-2 >2bc \Leftrightarrow b-2>bc \Leftrightarrow -2>b (c-1)$.
		    Since $b>0$ and $c\ge 1$, $b(c-1)\ge 0$; thus, this condition is never satisfied.
	    \end{compactitem}
      \end{compactitem}
  \item $\boxed{2} = \boxed{4}$ in Eq.~\eqref{eq:f1}:
      \begin{compactitem}
	  \item (i)-(ii):
	    \begin{align*}
		\lambda^c &= \frac{1}{l\lambda^{\frac{1}{b}}} \Leftrightarrow \lambda^{c+\frac{1}{b}=\frac{cb+1}{b}} = \frac{1}{l} \Leftrightarrow \lambda = \frac{1}{l^{\frac{b}{bc+1}}} \rightarrow 0,\text{ if } l\rightarrow \infty\quad [\Leftarrow \frac{b}{bc+1}>0].\\
		r(l,N) &= \frac{l^{\frac{3b}{bc+1}}\log^h(l)}{N^h} + \frac{1}{l^{\frac{bc}{bc+1}}} + \frac{l^{\frac{b}{bc+1}}}{l^2} \Rightarrow
		r(N) = \frac{\log^h(N)}{N^{h-\frac{3ab}{bc+1}}} + \frac{1}{N^{\frac{abc}{bc+1}}} + \frac{1}{N^{2a-\frac{ab}{bc+1}}}.
	    \end{align*}
	  Here, 
	  \begin{compactitem}
	      \item (ii):
		    $r(N)\rightarrow 0$, if 
		    \begin{compactitem}
			      \item $\ovalbox{1}\rightarrow 0$: Since $h>0$ we get $h-\frac{3ab}{bc+1}>0$, i.e., $\frac{h(bc+1)}{3b}>a$ using that $b>0$, $bc+1>0$.
			      \item $\ovalbox{2}\rightarrow 0$: $\frac{abc}{bc+1}>0$ -- the second condition is satisfied by our assumptions ($a>0$, $b>0$, $c>0$).
			      \item $\ovalbox{3}\rightarrow 0$: $2a-\frac{ab}{bc+1}>0$. Making use of the positivity of $a$ and $bc+1$, this requirement is equivalent to
			               $ 2 >\frac{b}{bc+1} \Leftrightarrow 2bc+2 > b \Leftrightarrow 2> b(1-2c)$,
				    which holds since $b(1-2c)<0$.
		    \end{compactitem}
		      Thus, we need $\frac{h(bc+1)}{3b} > a$. 
	      \item (i): $N^a  \left(\frac{1}{N^{\frac{ab}{bc+1}}}\right)^{\frac{b+1}{b}} \ge 1 \Leftrightarrow N^{a-\frac{a(b+1)}{bc+1}} \ge 1$.
		Thus it is sufficient to have $a-\frac{a(b+1)}{bc+1}>0 \Leftrightarrow 1 > \frac{b+1}{bc+1}$,
		using $a>0$.
	  \end{compactitem}
	    To sum up, for (i)-(ii) we got $\frac{h(bc+1)}{3b}>a$, $1 > \frac{b+1}{bc+1}$. 
	  \item (iii): 
	    \begin{compactitem}
		\item (i): $1 > \frac{b+1}{bc+1}$.
		\item $\ovalbox{2}\rightarrow 0$: no constraint.
		\item $\ovalbox{2} \ge \ovalbox{1}$: $\frac{1}{N^{\frac{abc}{bc+1}}} \ge \frac{\log^h(N)}{N^{h-\frac{3ab}{bc+1}}} \Leftrightarrow N^{h-\frac{3ab}{bc+1}-\frac{abc}{bc+1}} \ge \log^h(N)$.
		    Since $h>0$, this holds if $h-\frac{3ab}{bc+1}-\frac{abc}{bc+1} > 0 \Leftrightarrow h > a \frac{3b+bc}{bc+1} \Leftrightarrow \frac{h(bc+1)}{3b+bc} > a$,
		exploiting that $3b+bc>0$, $bc+1>0$.
		\item $\ovalbox{2} \ge \ovalbox{3}$: $\frac{1}{N^{\frac{abc}{bc+1}}} \ge \frac{1}{N^{2a-\frac{ab}{bc+1}}} \Leftrightarrow N^{2a-\frac{ab}{bc+1}-\frac{abc}{bc+1}}\ge 1$.
		    Hence, since $a>0$ and $bc+1>0$ we have
		    $2a-\frac{ab}{bc+1}-\frac{abc}{bc+1}>0\Leftrightarrow  2 > \frac{b+bc}{bc+1} \Leftrightarrow 2bc+2 > b +bc \Leftrightarrow bc+2 > b \Leftrightarrow 2> b(1-c)$.
		    This holds since $b(1-c)\le 0$.
	    \end{compactitem}
	    Thus, we got
	    \begin{align*}
		  \frac{h(bc+1)}{3b+bc} &> a, & 1 &> \frac{b+1}{bc+1} &
		r(N)&= \frac{1}{N^{\frac{abc}{bc+1}}} \rightarrow 0.
	    \end{align*}
      \end{compactitem}
  \item $\boxed{3} = \boxed{4}$ in Eq.~\eqref{eq:f1}:
      \begin{compactitem}
	  \item (i)-(ii):
	    \begin{align*}
		\frac{1}{l^2\lambda} &= \frac{1}{l\lambda^{\frac{1}{b}}} \Leftrightarrow \frac{1}{l} = \lambda^{1-\frac{1}{b}=\frac{b-1}{b}} \Leftrightarrow \frac{1}{l^{\frac{b}{b-1}}} = \lambda\rightarrow 0,\text{ if } l\rightarrow \infty\quad [\Leftarrow \frac{b}{b-1}>0].\\
		r(l,N) &= \frac{l^{\frac{3b}{b-1}}\log^h(l)}{N^h} + \frac{1}{l^{\frac{bc}{b-1}}} + \frac{l^{\frac{b}{b-1}}}{l^2} \Rightarrow
		r(N) = \frac{\log^h(N)}{N^{h-\frac{3ab}{b-1}}} + \frac{1}{N^{\frac{abc}{b-1}}} +  \frac{1}{N^{2a-\frac{ab}{b-1}}}.
	    \end{align*}
	    Here, 
	    \begin{compactitem}
		\item (ii):
		      $r(N)\rightarrow 0$ if 
		      \begin{compactitem}
			  \item $\ovalbox{1}\rightarrow 0$: Since $h>0$ we get $h-\frac{3ab}{b-1}>0$, i.e., $\frac{h(b-1)}{3b}>a$ using that $3b>0$ and $b-1>0$.
			  \item $\ovalbox{2}\rightarrow 0$: $\frac{abc}{b-1}>0$. This requirement holds by our assumptions [$a>0$, $b>1$, $c>0$].
			  \item $\ovalbox{3}\rightarrow 0$: $2a-\frac{ab}{b-1}>0$. By $a>0$ and $b-1>0$, this constraint is $2>\frac{b}{b-1} \Leftrightarrow 2b-2>b \Leftrightarrow b>2$.
		      \end{compactitem}
		      Hence, we need $\frac{h(b-1)}{3b}>a$, $b>2$.
	      \item (i): $N^a \left(\frac{1}{N^{\frac{ab}{b-1}}}\right)^{\frac{b+1}{b}} \ge 1 \Leftrightarrow N^{a-a\frac{b+1}{b-1}}  \ge 1$.
		Thus we need $a-a\frac{b+1}{b-1}>0\Leftrightarrow 1-1\frac{b+1}{b-1}>0 \Leftrightarrow 1>\frac{b+1}{b-1}$, where we used that $a>0$. 
		The $1>\frac{b+1}{b-1}$ is never satisfied since $\frac{b+1}{b-1}>1$.
	    \end{compactitem}
      \end{compactitem}
\end{compactitem}

\subsection{Numerical Experiments: Aerosol Prediction}\label{sec: numerical experiments}
In this section we provide numerical results to demonstrate the efficiency of the analysed ridge regression technique.
The experiments serve to illustrate that the MERR approach compares favourably to
\begin{enumerate}
    \item the only alternative, theoretically justified distribution regression method (since it avoids density estimation);\textsuperscript{\ref{footnote:RNDkitchensink}} see Section~\ref{sec:supHlearning},
    \item modern domain-specific, engineered methods (which beat state-of-the-art multiple instance learning alternatives); see Section~\ref{sec:aerosol prediction}.
\end{enumerate}
In our experiments  we used the ITE toolbox (Information Theoretical Estimators; \cite{szabo14information}).\footnote{The ITE toolbox contains the MERR method and its numerical demonstrations (among others); see \url{https://bitbucket.org/szzoli/ite/}.}

\subsubsection{Supervised entropy learning}\label{sec:supHlearning}
We compare  our MERR (RKHS based mean embedding ridge regression) algorithm
with  \cite{poczos13distribution}'s DFDR (kernel smoothing based distribution free distribution regression) method,
on a benchmark problem taken from the latter paper.
The goal is to learn the entropy of Gaussian distributions in a supervised way.
We chose an $A\in\R^{2\times 2}$ matrix, whose $A_{ij}$ entries were uniformly distributed on $[0,1]$ ($A_{ij}\sim U[0,1]$).
We constructed $100$ sample sets from $\{N(0,\Sigma_u)\}_{u=1}^{100}$, where $\Sigma_u=R(\beta_u)AA^TR(\beta_u)^T$ and $R(\beta_u)$ was a 2d rotation matrix with angle $\beta_u\sim U[0,\pi]$.
\begin{figure}
\centering
\subfloat[][Entropy of Gaussian]{\includegraphics[width = 0.35\linewidth]{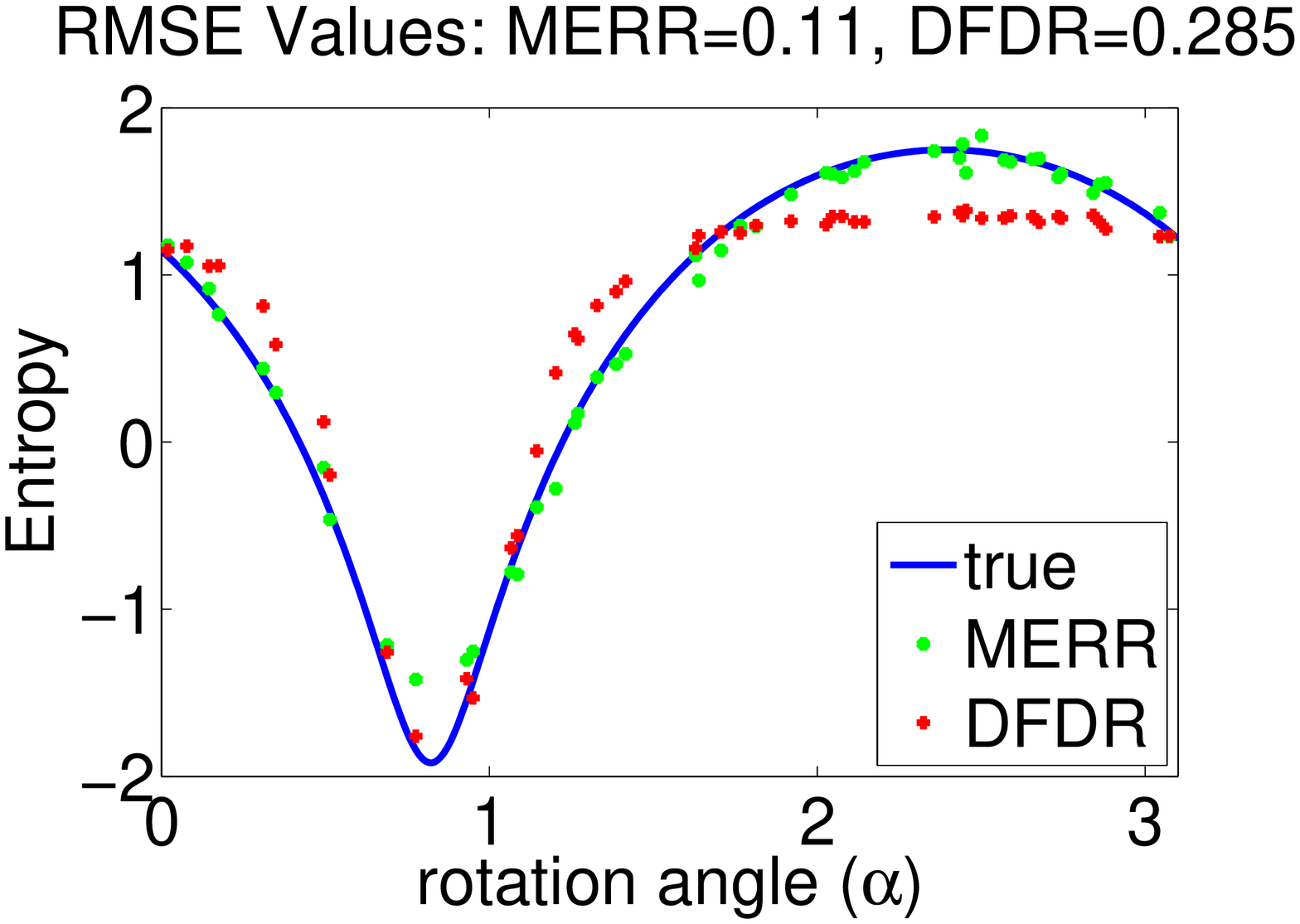}}
\subfloat[][Boxplot of RMSE]{\includegraphics[trim= 0 26 0 0, width = 0.35\linewidth]{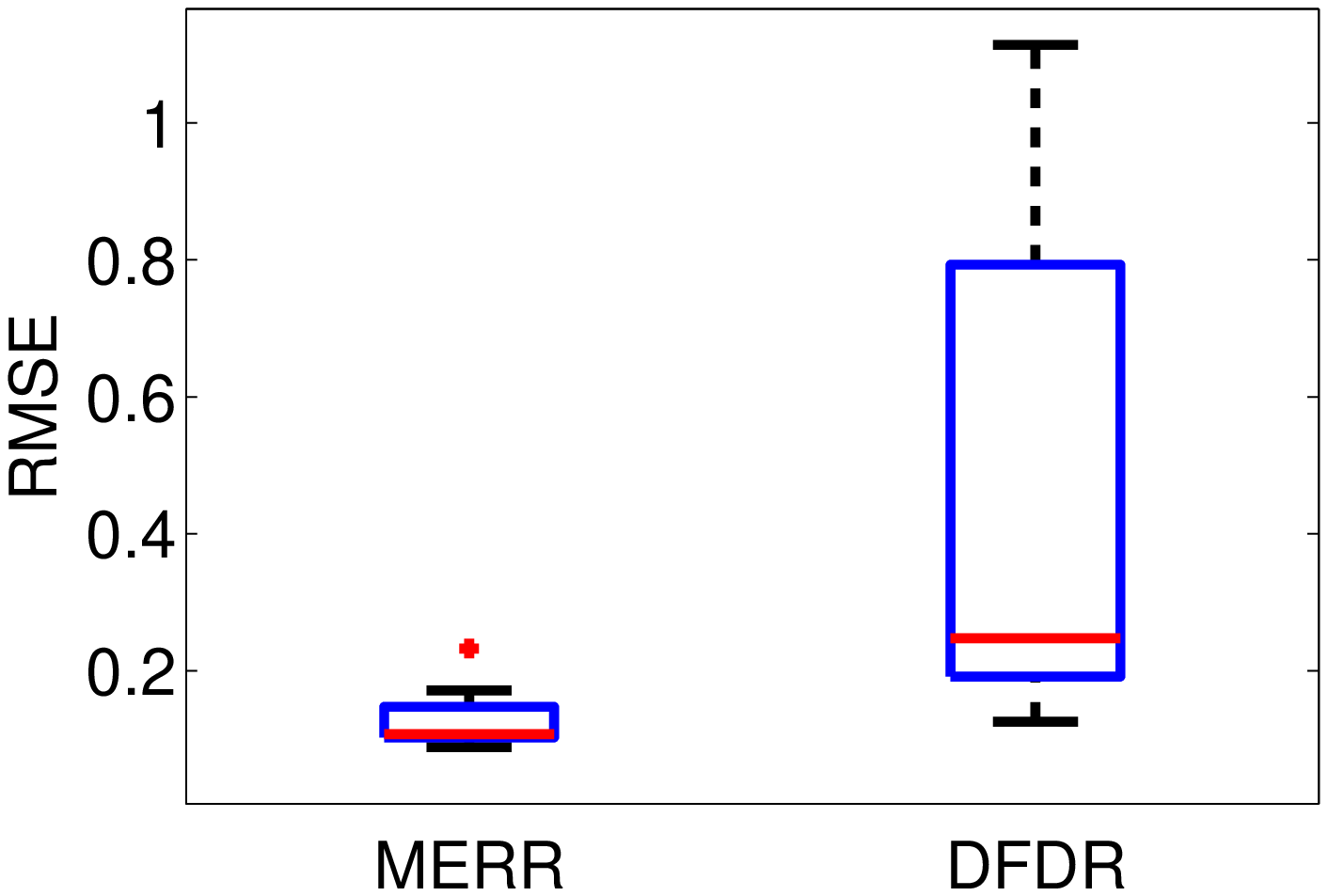}}
\caption{(a) Learned entropy of a one-dimensional marginal distribution of a rotated 2d Gaussian. Axes $x$: rotation angle in $[0,\pi]$. Axis $y$: entropy.
(b) RMSE values of the MERR and DFDR algorithms. Boxplots are calculated from $25$ experiments. \label{fig:boxplot_RMSE}}
\end{figure}
From each $N(0,\Sigma_u)$ distribution we sampled $500$ 2-dimensional i.i.d.\ points. From the  $100$ sample sets, $25$  were used for training, $25$ for validation
(i.e., selecting appropriate parameters), and $50$ points for testing. Our goal is to learn the entropy of the first marginal distribution: $H=\frac{1}{2}\ln(2\pi e\sigma^2)$,
where $\sigma^2=M_{1,1}$, $M=\Sigma_u\in \R^{2 \times 2}$. Figure~\ref{fig:boxplot_RMSE}(a) displays the learned entropies of the $50$ test sample sets in a
typical experiment. We compare the results of DFDR and MERR.
One can see that the true and the estimated values are close to each other for both algorithms, but MERR performs better.
The boxplot diagrams of the RMSE (root mean square error) values calculated from $25$ experiments confirm this performance advantage (Figure~\ref{fig:boxplot_RMSE}(b)).
A reason why MERR achieves better performance is that DFDR needs to do many density estimations, which can be very challenging when the sample sizes are small.
By contrast, the MERR algorithm  does not require density estimation.

\subsubsection{Aerosol prediction} \label{sec:aerosol prediction}
Aerosol prediction is one of the largest challenges of current climate research; we chose this problem as a further testbed of our method.
 \cite{wang12mixture} pose the AOD (aerosol optical depth) prediction problem as a MIL task: (i)
a given pixel of a multispectral image corresponds to a small area of $200\times 200m^2$, (ii)
spatial variability of AOD can be considered to be small over distances up to $100km$,
(iii) ground-based instruments provide AOD labels ($y_i\in\R$), (iv) a bag consists of randomly selected pixels within a
$20km$ radius around an AOD sensor. The MIL task can be tackled using our MERR approach, assuming that (i) bags correspond to distributions ($x_i$),
(ii) instances in the bag ($\{x_{i,n}\}_{n=1}^N$) are samples from the distribution.

We selected the MISR1 dataset \cite{wang12mixture}, where (i) cloudy pixels are also included, (ii) there are $800$ bags with (iii) $100$ instances in each bag, (iv) the
instances are 16-dimensional ($x_{i,n}\in\R^{16}$). Our baselines are the reported state-of-the-art EM (expectation-maximization) methods achieving average $100\times RMSE = 7.5-8.5$ ($\pm 0.1-0.6$) accuracy. The experimental protocol followed the original work, where 5-fold cross-validation ($4\times 160$ ($160$) samples for
training (testing)) was repeated $10$ times; the only difference is that we made the problem a bit harder, as we used only $3\times 160$ samples for training,
$160$ for validation (i.e., setting the $\lambda$ regularization and the $\theta$ kernel parameter), and $160$ for testing.

\begin{itemize}
  \item Linear $K$: In the first set of experiments, $K$ was linear. To study the robustness of our method, we picked $10$ different kernels ($k$) and their ensembles: the Gaussian, exponential, Cauchy, generalized t-student,
	  polynomial kernel of order $2$ and $3$ ($p=2$ and $3$), rational quadratic, inverse multiquadratic kernel, Mat{\'e}rn kernel (with
	  $\frac{3}{2}$ and $\frac{5}{2}$ smoothness parameters). The expressions for these kernels are
	  \begin{align*}
	    k_G(a,b) &= e^{-\frac{\left\|a-b\right\|_2^2}{2\theta^2}},       & k_e(a,b) &= e^{-\frac{\left\|a-b\right\|_2}{2\theta^2}}, &  k_C(a,b) &= \frac{1}{1+\frac{\left\|a-b\right\|_2^2}{\theta^2}},\\
	    k_t(a,b) &= \frac{1}{1+\left\|a-b\right\|_2^{\theta}},           & k_p(a,b) &= \left(\left<a,b\right>+\theta\right)^p,      & k_r(a,b) &= 1-\frac{\left\|a-b\right\|_2^2}{\left\|a-b\right\|_2^2+\theta},\\
	    k_i(a,b) &= \frac{1}{\sqrt{\left\|a-b\right\|_2^2+\theta^2}},    &k_{M,\frac{3}{2}}(a,b) &= \left(1+\frac{\sqrt{3}\left\|a-b\right\|_2}{\theta}\right)e^{-\frac{\sqrt{3}\left\|a-b\right\|_2}{\theta}}&\\
	  \end{align*}
	  \vspace*{-1.2cm}
	  \begin{align*}
	  \hspace*{-4.7cm} k_{M,\frac{5}{2}}(a,b) &= \left(1+\frac{\sqrt{5}\left\|a-b\right\|_2}{\theta}+\frac{5\left\|a-b\right\|_2^2}{3\theta^2}\right)e^{-\frac{\sqrt{5}\left\|a-b\right\|_2}{\theta}},
	  \end{align*}
	  where $p=2,3$ and $\theta>0$. The explored parameter domain was $(\lambda, \theta)\in \{2^{-65},2^{-64},\ldots,2^{-3}\}\times \{2^{-15},2^{-14},\ldots,2^{10}\}$; increasing the domain further did not improve the results.

	  Our results are summarized in Table~\ref{tab:MISR1}. According to the table,
	  we achieve $100\times RMSE = 7.91$ ($\pm 1.61$) using a single kernel, or
	  $7.86$ ($\pm 1.71$) with ensemble of kernels (further performance improvements might be obtained by learning the weights).
	  \begin{table}
	    %\small
	    \begin{center}
	    \caption{Prediction accuracy of the MERR method in AOD prediction using different kernels: $100\times RMSE (\pm std)$. $K$: linear. The best single and ensemble predictions are written in bold.}\label{tab:MISR1}
	    \begin{tabular}{@{}l@{\hspace{0.25cm}}l@{\hspace{0.25cm}}l@{\hspace{0.25cm}}l@{\hspace{0.25cm}}l@{\hspace{0.25cm}}l@{}}
	      \toprule
	      $k_G$               & $k_e$                  & $k_C$               & $k_t$ &     $k_p (p=2)$         & $k_p (p=3)$                            \\
	      $7.97$ ($\pm 1.81$) & $8.25$ ($\pm 1.92$)    &$7.92$ ($\pm 1.69$)  & $8.73$ ($\pm 2.18$) & $12.5$ ($\pm 2.63$) & $171.24$ ($\pm 56.66$) \\\midrule
	      $k_r$ & $k_i$ & $k_{M,\frac{3}{2}}$ & $k_{M,\frac{5}{2}}$    &ensemble & \\
	      $9.66$ ($\pm 2.68$) & $\tb{7.91}$ ($\pm \tb{1.61})$ &$8.05$ ($\pm 1.83$) & $7.98$ ($\pm 1.75$)    & $\tb{7.86}$ ($\pm \tb{1.71}$)\\
	      \bottomrule
	    \end{tabular}
	    \end{center}
	  \end{table}
  \item Nonlinear $K$: We also studied the efficiency of nonlinear $K$-s. In this case, the argument of $K$ was $\left\|\mu_a-\mu_b\right\|_H$ instead of $\left\|a-b\right\|_2$ (see the definition of $k$-s); for $K$ examples, see Table~\ref{tab:K examples}. Our obtained results are
	summarized in Table~\ref{tab:log2:nonlinK}. One can see that using nonlinear $K$ kernels, the RMSE error drops to $7.90$  ($\pm 1.63$) in the single prediction case, and decreases further to
	$7.81$ ($\pm 1.64$) in the ensemble setting.
      \begin{table}
      \begin{center}
      \caption{Prediction accuracy of the MERR method in AOD prediction using different kernels: $100\times RMSE (\pm std)$; single prediction case. $K$: nonlinear. Rows: kernel $k$. Columns: kernel $K$. For each row ($k$), the smallest RMSE value is written in bold.} \label{tab:log2:nonlinK}
      \begin{tabular}{@{}lccccc@{}}\toprule
      $k\backslash K$        &$K_G$                 &$K_e$               &$K_C$                &$K_t$                                 &$K_{M,\frac{3}{2}}$                     \\\hline
      $k_e$                  &$8.14$ ($\pm 1.80$)   &$8.10$ ($\pm 1.81$) &$8.14$ ($\pm 1.81$)  &$\mathbf{8.07}$ ($\pm \mathbf{1.77}$) &$8.09$           ($\pm 1.88$)           \\
      $k_C$                  &$7.97$ ($\pm 1.58$)   &$8.13$ ($\pm 1.79$) &$7.96$ ($\pm 1.62$)  &$8.09$          ($\pm 1.69$)          &$\mathbf{7.90}$  ($\pm \mathbf{1.63}$)  \\
      $k_{M,\frac{3}{2}}$    &$8.00$ ($\pm 1.66$)   &$8.14$ ($\pm 1.80$) &$8.00$ ($\pm 1.69$)  &$8.08$          ($\pm 1.72$)          &$\mathbf{7.96}$  ($\pm  \mathbf{1.69}$) \\
      $k_i$                  &$8.01$ ($\pm 1.53$)   &$8.17$ ($\pm 1.74$) &$8.03$ ($\pm 1.63$)  &$7.93$          ($\pm 1.57$)          &$8.04$           ($\pm 1.67$)           \\\midrule
      $k\backslash K$        &$K_{M,\frac{5}{2}}$   &$K_r$               &$K_i$                & linear \\\hline
      $k_e$                  &$8.14$ ($\pm 1.78$)   &$8.12$ ($\pm 1.81$) &$8.12$ ($\pm 1.80$)  & $8.25$ ($\pm{1.92}$)\\
      $k_C$                  &$7.95$ ($\pm 1.60$)   &$7.92$ ($\pm 1.61$) &$7.93$ ($\pm 1.61$)  & $7.92$ ($\pm 1.69$)\\
      $k_{M,\frac{3}{2}}$    &$8.02$ ($\pm 1.71$)   &$8.04$ ($\pm 1.69$) &$7.98$ ($\pm 1.72$)  & $8.05$ $(\pm 1.83)$\\
      $k_i$                  &$8.05$ ($\pm 1.61$)   &$8.05$ ($\pm 1.63$) &$8.06$ ($\pm 1.65$)  & $\mathbf{7.91} (\pm \mathbf{1.61})$\\\bottomrule
      \end{tabular}
      \end{center}
      \end{table}
\end{itemize}

Despite the fact that MERR has no domain-specific knowledge wired in, the results fall within the same range as \cite{wang12mixture}'s algorithms.
The prediction is fairly precise and robust to the choice of the kernel, however polynomial kernels perform poorly (they violate our boundedness assumption).

\end{document}